\documentclass{emsprocart}
\usepackage{xypic}

\contact[ebeling@math.uni-hannover.de]{Wolfgang Ebeling, Institut f\"{u}r Algebraische Geometrie, Leibniz Universit\"{a}t Hannover, 
Postfach 6009, D-30060 Hannover, Germany}

\newtheorem{theorem}{Theorem}[section]
\newtheorem{corollary}[theorem]{Corollary}

\newtheorem{proposition}[theorem]{Proposition}
\newtheorem{conjecture}[theorem]{Conjecture}

\theoremstyle{definition}
\newtheorem{definition}[theorem]{Definition}

\newcommand{\CC}{{\mathbb C}}
\newcommand{\HH}{{\mathbb H}}
\newcommand{\PP}{{\mathbb P}}
\newcommand{\QQ}{{\mathbb Q}}

\newcommand{\ZZ}{{\mathbb Z}}

\newcommand{\calO}{{\cal O}}

\newcommand{\calC}{{\cal C}}

\newcommand{\eps}{\varepsilon}

\def\H{{\mathcal H}}
\newcommand{\Conjsub}{{\rm Conjsub}\,}
\newcommand{\Ext}{\mathrm{Ext}}

\newcommand{\TTT}{\mathsf{T}}
\newcommand{\tytt}[1]{{\!#1}}  
\newcommand{\ff}{{\bf f}}
\newcommand{\hh}{{\bf h}}

\title[Homological mirror symmetry for singularities]{Homological mirror symmetry for singularities}

\author[Wolfgang Ebeling]{Wolfgang Ebeling\thanks{Partially supported by the DFG-programme SPP1388 ''Representation Theory''}}

\begin{document}

\begin{abstract} We give a survey on results related to the Berglund-H\"ubsch duality of invertible polynomials and the homological mirror symmetry conjecture for singularities.
\end{abstract}



\maketitle

\section{Introduction}
V.~I.~Arnold observed a strange duality between the 14 exceptional unimodal singularities. When physicists came up with the idea of mirror symmetry, it was found that Arnold's strange duality can be considered as part of the mirror symmetry of K3 surfaces. In his 1994 talk at the International Congress of Mathematicians, M.~Kontsevich \cite{Kon} proposed an interpretation of the mirror phenomenon in mathematical terms which is now commonly known as the {\em homological mirror symmetry conjecture}. It was originally formulated for two mirror symmetric Calabi-Yau manifolds $X$ and $X'$ and states that there is an equivalence between the derived category of coherent sheaves on $X$ and the derived Fukaya category of $X'$ and vice versa. 

Kontsevich also suggested that homological mirror symmetry can be extended to a more general setting by considering Landau-Ginzburg models. Many aspects of Landau-Ginzburg models are related to singularity theory. One of the early constructions of mirror symmetric manifolds was the construction of P.~Berglund and T.~H\"ubsch \cite{BH}. They considered a polynomial $f$ of a special form, a so called {\em invertible}Ê one, and its {\em Berglund-H\"ubsch transpose} $\widetilde{f}$: see Sect.~\ref{sect:inv}. These polynomials can be considered as potentials of Landau-Ginzburg models. This construction can also be generalized to an orbifold setting.  One can formulate different versions of the homological mirror symmetry conjecture for Berglund-H\"ubsch pairs. 
It turned out that Arnold's strange duality is also part of this duality and features of Arnold's strange duality appeared as features of homological mirror symmetry.
We review results related to these conjectures.

We briefly outline the contents of this survey. We start by discussing Arnold's strange duality. In Sect.~\ref{sect:inv}, we review the notion of an invertible polynomial and the  Berglund-H\"ubsch construction. In Sect.~\ref{sect:HMS}, we state the homological mirror symmetry conjectures for invertible polynomials. In Sect.~\ref{sect:SD}, we give a survey on the evidence for these conjectures. More precisely, we give a generalization of Arnold's strange duality. In Sect.~\ref{sect:E}, we show that the mirror symmetry for Berglund-H\"ubsch dual pairs also holds on the level of suitably defined Hodge numbers. For this purpose we discuss the notion of an orbifold E-function of a polynomial with an isolated singularity at the origin and we consider these functions for dual pairs. Another feature of Arnold's strange duality was discovered by K.~Saito and is known as Saito duality. We discuss how this duality generalizes to the Berglund-H\"ubsch duality. In Sect.~\ref{sect:Ex}, we compile the more detailed information one has about specific classes of singularities, like the simple, unimodal and bimodal singularities. Finally, we derive in Sect.~\ref{sect:EW} the extension of Arnold's strange duality involving complete intersection singularities \cite{EW} from the Berglund-H\"ubsch construction.

\section{Arnold's strange duality} \label{sect:Arnold}
According to Arnold's classification of singularities \cite{Ar}, there are 14 exceptional unimodal singularities. Setting the modulus equal to zero, they can be given by equations $f(x,y,z)=0$ where the polynomial $f$ is given in Table~\ref{TabArnold}. We use the name of Arnold for the corresponding singularity.

We associate Dolgachev and Gabrielov numbers to these singularities as follows. 

Consider the quotient stack
\[
{\mathcal C}_{f}:=\left[f^{-1}(0)\backslash\{0\}\left/\CC^\ast \right.\right].
\]
This is a Deligne-Mumford stack and can be regarded as a smooth projective line $\PP^1$ with  three isotropic points of orders $\alpha_1, \alpha_2, \alpha_3$. The numbers $(\alpha_1, \alpha_2, \alpha_3)$ are called the {\em Dolgachev numbers} of $f$ \cite{Dolgachev74, Dolgachev75}.

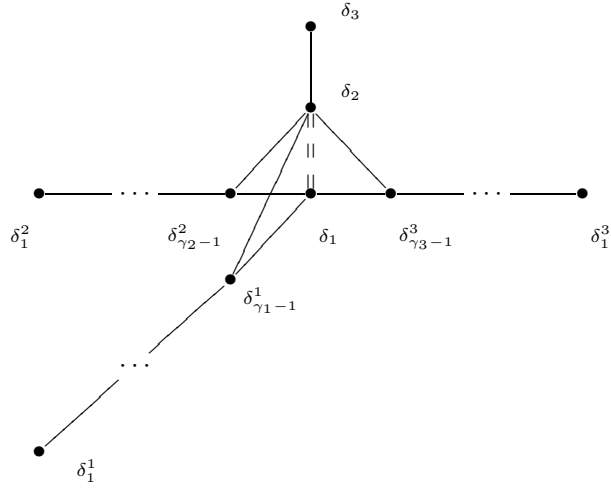
\begin{figure}
$$
\xymatrix{ 
 & & & *{\bullet} \ar@{-}[d] \ar@{}^{\delta_{3}}[r] & & & \\
 & & & *{\bullet} \ar@{==}[d] \ar@{-}[dr]  \ar@{-}[ldd] \ar@{}^{\delta_{2}}[r]
 & & &  \\
 *{\bullet} \ar@{-}[r] \ar@{}_{\delta^2_1}[d]  & {\cdots} \ar@{-}[r]  & *{\bullet} \ar@{-}[r] \ar@{-}[ur]   \ar@{}_{\delta^2_{\gamma_2-1}}[d] & *{\bullet} \ar@{-}[dl] \ar@{-}[r] \ar@{}^{\delta_{1}}[d] & *{\bullet} \ar@{-}[r]  \ar@{}^{\delta^3_{\gamma_3-1}}[d]  & {\cdots} \ar@{-}[r]  &*{\bullet} \ar@{}^{\delta^3_1}[d]   \\
 & &   *{\bullet} \ar@{-}[dl] \ar@{}_{\delta^1_{\gamma_1-1}}[r]  & & & & \\
 & {\cdots} \ar@{-}[dl] & & & & & \\
*{\bullet}  \ar@{}_{\delta^1_1}[r] & & & & & &
  } 
$$
\caption{The graph $S_{\gamma_1,\gamma_2,\gamma_3}$} \label{FigSpqr}
\end{figure}

The manifold $V_f:=f^{-1}(1)$ is called the {\em Milnor fibre} of $f$. Since $f$ has an isolated singularity at the origin, the only interesting homology group is $H_2(V_f,\ZZ)$. We denote by $\langle \ , \ \rangle$ the intersection form on $H_2(V_f,\ZZ)$
and by $H=(H_2(V_f,\ZZ),\langle \ , \ \rangle)$ the Milnor lattice.
A.~M.~Gabrielov
\cite{Gabrielov74} has shown that there exists a weakly distinguished basis of
vanishing cycles of
$H$ with a Coxeter-Dynkin diagram of the form of Fig.~\ref{FigSpqr}. The author
\cite{Ebeling81} (see also \cite{Ebeling96}) has shown that one can even find a distinguished basis
\[(\delta_1, \delta^1_1, \ldots \delta^1_{\gamma_1-1}, \delta^2_1, \ldots, \delta^2_{\gamma_2-1}, \delta^3_1, \ldots, \delta^3_{\gamma_3-1}, \delta_2, \delta_3)\]
with this Coxeter-Dynkin diagram. (For the notions of a
distinguished and weakly distinguished basis of vanishing cycles see e.g.\
\cite{Eb}). The numbers $\gamma_1$, $\gamma_2$, $\gamma_3$ are called the {\em Gabrielov
numbers} of the singularity. Here each vertex represents a
sphere of self-intersection number $-2$, two vertices connected by a
single solid edge have intersection number 1, two vertices connected by a
double broken line have intersection number $-2$ and vertices which are not connected have intersection number 0. 

Arnold \cite{Ar} has now observed: There exists an involution $X \mapsto X^\vee$ (indicated in Table~\ref{TabArnold}) on the
set of the 14 exceptional unimodal singularities, such that the Dolgachev numbers of $X$ are the Gabrielov numbers of $X^\vee$ and the Gabrielov numbers of $X$ are the Dolgachev numbers of $X^\vee$. This is called {\em Arnold's strange duality}. 

\begin{table}
\begin{center}
\begin{tabular}{ccccc}  \hline
Name  & $f$ & $\alpha_1,\alpha_2, \alpha_3$ & $\gamma_1, \gamma_2, \gamma_3$ & Dual
\\ \hline
$E_{12}$  & $x^2+y^3+z^7$ & $2,3,7$ & $2,3,7$ & $E_{12}$ \\ 
$E_{13}$  & $x^2+y^3+yz^5$ & $2,4,5$ & $2,3,8$ & $Z_{11}$ \\ 
$E_{14}$   & $x^3+y^2+yz^4$ & $3,3,4$ & $2,3,9$ & $Q_{10}$ \\ 
$Z_{11}$ & $x^2+zy^3+z^5$ & $2,3,8$ & $2,4,5$ & $E_{13}$ \\
$Z_{12}$  & $x^2+zy^3+yz^4$ & $2,4,6$ & $2,4,6$ & $Z_{12}$ \\
$Z_{13}$ & $x^2+xy^3+yz^3$ & $3,3,5$ & $2,4,7$ & $Q_{11}$ \\ 
$Q_{10}$ &  $x^3+zy^2+z^4$ & $2,3,9$ & $3,3,4$ & $E_{14}$ \\
$Q_{11}$ & $x^2y+y^3z+z^3$ & $2,4,7$ & $3,3,5$ & $Z_{13}$ \\ 
$Q_{12}$  & $x^3+zy^2+yz^3$ & $3,3,6$ & $3,3,6$ &  $Q_{12}$ \\ 
$W_{12}$  & $x^5 + y^2+yz^2$ & $2,5,5$ & $2,5,5$ & $W_{12}$ \\ 
$W_{13}$  & $x^2+xy^2+yz^4$ & $3,4,4$ & $2,5,6$ & $S_{11}$ \\
$S_{11}$ & $x^2y+y^2z+z^4$ & $2,5,6$ & $3,4,4$ & $W_{13}$ \\
$S_{12}$  & $x^3y+y^2z+z^2x$ & $3,4,5$ & $3,4,5$ &  $S_{12}$ \\ 
$U_{12}$  &  $x^4+zy^2+yz^2$ & $4,4,4$ & $4,4,4$ & $U_{12}$ \\ 
\hline
\end{tabular}
\end{center}
\caption{Arnold's strange duality}\label{TabArnold}
\end{table}

Consider $f$ as a function $f: (\CC^3,0) \to (\CC,0)$. A characteristic homeomorphism of
the Milnor fibration of $f$ induces an automorphism $c:H_2(V_f,\ZZ) \to
H_2(V_f,\ZZ)$ called the {\em (classical) monodromy operator} of $f$. It is the {\em Coxeter element} corresponding to a distinguished basis $\{ \delta_1, \ldots , \delta_\mu\}$ of vanishing cycles of $f$. By this we mean the following: Each vanishing cycle $\delta_i$ defines a reflection
\[\begin{array}{cccc} s_{\delta_i}: &  H_2(V_f,\ZZ) & \to & H_2(V_f,\ZZ)\\
 & x & \mapsto & s_{\delta_i}(x):= x - \frac{2\langle x,\delta_i \rangle}{\langle \delta_i, \delta_i \rangle} \delta_i \end{array}
\]
Then 
\[ c= s_{\delta_1} \circ s_{\delta_2} \circ \cdots \circ s_{\delta_\mu}.
\]
It is a well known theorem  (see e.g. \cite{Brieskorn70}) that the eigenvalues of
$c$ are roots of unity. This means that the characteristic polynomial
$\phi(\lambda) = \det (\lambda I - c)$ of $c$ is a monic polynomial the roots of
which are roots of unity. Moreover, since $f$ is weighted homogeneous, the operator $c$ has finite order $h$. Such a polynomial can be written in the form
$$\phi(\lambda)= \prod_{m | h} (\lambda^m -1)^{\chi_m} \quad \mbox{for} \
\chi_m \in \ZZ.$$
K.~Saito
\cite{Sa1, Sa2} defines a {\em dual polynomial}
$\phi^\vee(\lambda)$ to $\phi(\lambda)$:
$$\phi^\vee(\lambda) = \prod_{k | h} (\lambda^k -1)^{-\chi_{h/k}}.$$
He obtains the following result.

\begin{theorem}[Saito] \label{thm:Saito}
If $\phi(\lambda)$ is the characteristic polynomial of the monodromy of an
 exceptional unimodal singularity $X$, then $\phi^\vee(\lambda)$ is the
corresponding polynomial of the dual singularity $X^\vee$.
\end{theorem}

The author and C.T.C.~Wall \cite{EW} discovered an extension of Arnold's strange duality embracing on one hand series of bimodal singularities and on the other, isolated complete intersection singularities (ICIS) in $\CC^4$. The duals of the complete intersection singularities are not themselves singularities, but are virtual ($k=-1$) cases of series (e.g.\  $W_{1,k} : k \geq 0$) of bimodal singularities. They associated to these Dolgachev and Gabrielov numbers and showed that all the features of Arnold's strange duality continue to hold. Moreover, in \cite{Ebeling99} the author showed that also Saito's duality holds for this duality. We come back to this extension in Sect.~\ref{sect:EW}.

\section{Invertible polynomials} \label{sect:inv}
We recall some general definitions about invertible polynomials.

Let $f(x_1,\dots, x_n)$ be a  weighted homogeneous polynomial, namely, a polynomial with the property that there are positive integers $w_1,\dots ,w_n$ and $d$ such that 
$f(\lambda^{w_1} x_1, \dots, \lambda^{w_n} x_n) = \lambda^d f(x_1,\dots ,x_n)$ 
for $\lambda \in \CC^\ast$. We call $(w_1,\dots ,w_n;d)$ a system of {\em weights}. 
\begin{definition}
A  weighted homogeneous polynomial $f(x_1,\dots ,x_n)$ is called {\em invertible} if 
the following conditions are satisfied:
\begin{enumerate}
\item the number of variables ($=n$) coincides with the number of monomials 
in the polynomial $f(x_1,\dots x_n)$, 
namely, 
\[
f(x_1,\dots ,x_n)=\sum_{i=1}^na_i\prod_{j=1}^nx_j^{E_{ij}}
\]
for some coefficients $a_i\in\CC^\ast$ and non-negative integers 
$E_{ij}$ for $i,j=1,\dots, n$,
\item the system of weights $(w_1,\dots ,w_n;d)$ of $f$ is uniquely determined by 
the polynomial $f(x_1,\dots ,x_n)$ up to a constant factor ${\rm gcd}(w_1,\dots ,w_n;d)$, 
namely, the matrix $E:=(E_{ij})$ is invertible over $\QQ$.
\end{enumerate}
An invertible polynomial is called {\em non-degenerate}, if it has an isolated singularity at the origin. 
\end{definition}

Without loss of generality we shall assume that $\det E > 0$.

An invertible polynomial has a {\em canonical system of weights} $W_f=(w_1, \ldots , w_n;d)$ given by the unique solution of the equation
\begin{equation*}
E
\begin{pmatrix}
w_1\\
\vdots\\
w_n
\end{pmatrix}
={\rm det}(E)
\begin{pmatrix}
1\\
\vdots\\
1
\end{pmatrix}
,\quad 
d:={\rm det}(E).
\end{equation*}
This system of weights is in general non-reduced, i.e.\  in general 
\[ c_f:= {\rm gcd}(w_1, \ldots , w_n,d)>1.
\]

\begin{definition} Let $h(x_1, \ldots, x_n)$ be any polynomial. Let $G_h$ be the (finite) group of diagonal symmetries of $h$, i.e. 
\[
G_h := \left\{ (\lambda_1, \ldots , \lambda_n) )\in(\CC^\ast)^n \, \left| \, h(\lambda_1x_1, \ldots, \lambda_n x_n) = h(x_1, \ldots, x_n) \right\} \right..
\]
\end{definition}

\begin{definition}
Let $f(x_1,\dots ,x_n)=\sum_{i=1}^na_i\prod_{j=1}^nx_j^{E_{ij}}$ be an invertible polynomial. Consider the free abelian group $\oplus_{i=1}^n\ZZ\vec{x_i}\oplus \ZZ\vec{f}$ 
generated by the symbols $\vec{x_i}$ for the variables $x_i$ for $i=1,\dots, n$
and the symbol $\vec{f}$ for the polynomial $f$.
The {\em maximal grading} $L_f$ of the invertible polynomial $f$ 
is the abelian group defined by the quotient 
\[
L_f:=\bigoplus_{i=1}^n\ZZ\vec{x_i}\oplus \ZZ\vec{f}\left/I_f\right.,
\]
where $I_f$ is the subgroup generated by the elements 
\[
\vec{f}-\sum_{j=1}^nE_{ij}\vec{x_j},\quad i=1,\dots ,n.
\]
\end{definition}

\begin{definition}
Let $f(x_1,\dots ,x_n)$ be an invertible polynomial and $L_f$ be the maximal grading of $f$.
The {\em maximal abelian symmetry group} $\widehat{G}_f$ of $f$ is the abelian group defined by 
\[
\widehat{G}_f:={\rm Spec}(\CC L_f),
\]
where $\CC L_f$ denotes the group ring of $L_f$. Equivalently, 
\[
\widehat{G}_f=\left\{(\lambda_1,\dots ,\lambda_n)\in(\CC^\ast)^n \, \left| \,
\prod_{j=1}^n \lambda_j ^{E_{1j}}=\dots =\prod_{j=1}^n\lambda_j^{E_{nj}}\right\} \right..
\]
\end{definition}

We have
\[
G_f=\left\{(\lambda_1,\dots ,\lambda_n)\in \widehat{G}_f \, \left| \, \prod_{j=1}^n \lambda_j ^{E_{1j}}=\dots =\prod_{j=1}^n\lambda_j^{E_{nj}}=1 \right\} \right..
\]

Let $f(x_1,\dots ,x_n)$ be an invertible polynomial and $W_f=(w_1, \ldots , w_n;d)$ be the canonical system of weights associated to $f$. Set 
\[ q_i := \frac{w_i}{d}, \quad i=1, \ldots , n.
\]
Note that $G_f$ always contains the {\em exponential grading operator}
\[
g_0:=({\bf e}[q_1], \ldots , {\bf e}[q_n]),
\]
where we use the notation ${\bf e}[-] := \exp(2\pi \sqrt{-1} \cdot  -)$.
Let $G_0$ be the subgroup of $G_f$ generated by $g_0$. One has (cf.\ \cite{ET2})
\[
[G_f : G_0]=c_f.
\]
\begin{definition}
Let $f(x_1, \ldots , x_n)$ be a weighted homogeneous polynomial with reduced system of weights
$W=(w_1, \ldots ,w_n;d)$. 
The integer
\[
a_f := d- \sum_{i=1}^n w_i
\]
is called the {\em Gorenstein parameter} of $f$. It is also usual to denote $\epsilon_f:= -a_f$ the Gorenstein parameter of $f$, see e.g.\ \cite{Sa2}. 
\end{definition}

Let $f(x_1,\dots ,x_n)=\sum_{i=1}^na_i\prod_{j=1}^nx_j^{E_{ij}}$ be an invertible polynomial. 
\begin{definition}[Berglund, H\"ubsch]
Following \cite{BH}, the {\em Berglund-H\"ubsch transpose} of $\widetilde{f}(x_1, \ldots , x_n)$ of $f$ is defined by 
\[ 
\widetilde{f}(x_1,\dots ,x_n)=\sum_{i=1}^na_i\prod_{j=1}^nx_j^{E_{ji}}.
\]
\end{definition}

\begin{definition}[Berglund, Henningson]
By \cite{BHe}, for a subgroup $G \subset G_f$ its {\em dual group}  $\widetilde{G}$ is defined by 
\[
\widetilde{G}:= {\rm Hom}(G_f/G, \CC^\ast).
\]
\end{definition}

One has the following easy facts:
\begin{itemize}
\item $\widetilde{G_f} = \{ e \}$
\item $H \subset G \Rightarrow \widetilde{G} \subset \widetilde{H}$
\item $\widetilde{\widetilde{H}} = H$
\end{itemize}

Note that ${\rm Hom}(G_f/G, \CC^\ast)$ is isomorphic to $G_{\widetilde{f}}$, see \cite{BHe}.
By \cite{Kr} (see also \cite[Lemma 1]{EGPEMS}), we have
\[
\widetilde{G}_0 = {\rm SL}_n(\ZZ) \cap G_{\widetilde{f}}. 
\]
Moreover, by \cite[Proposition 3.1]{ET2}, we have $|\widetilde{G}_0|=c_f$.

For a subgroup  $G \subset G_f$, let $\widehat{G}$ be the subgroup of $\widehat{G}_f$ 
defined by the following commutative diagram of short exact sequences
\[ 
\xymatrix{
\{ 1 \}\ar[r] & G \ar[r]\ar@{^{(}->}[d]  & \widehat{G} \ar[r]\ar@{^{(}->}[d] &  \CC^\ast \ar[r]\ar@{=}[d]& \{ 1 \}\\
\{ 1 \}\ar[r] & G_f \ar[r] & \widehat{G}_f \ar[r] & \CC^\ast \ar[r]& \{ 1 \}
}.
\]

\section{Homological mirror symmetry} \label{sect:HMS}
There are several versions of the homological mirror symmetry conjecture for singularities. 

Let $f(x,y,z)$ be a polynomial which has an isolated singularity at the origin. A distinguished basis of vanishing cycles
in the Milnor fiber of $f$ can be categorified to an $A_\infty$-category ${\rm Fuk}^{\to}(f)$
called the directed Fukaya category.  Any two distinguished bases of vanishing cycles are connected by a sequence of {\em Gabrielov transformations} \cite{Gabrielov73}. The set of objects of ${\rm Fuk}^{\to}(f)$ is a distinguished basis of (Lagrangian) vanishing cycles and the spaces of morphisms are Lagrangian intersection Floer complexes. It can be shown that Gabrielov transformations correspond to {\em mutations} of the category (\cite{Se1}, see also e.g.\ \cite{KY}). Since different choices of distinguished bases are related by mutations, the derived category $D^b{\rm Fuk}^\to(f)$ is independent of this choice and 
is therefore, as a triangulated category, an invariant of the polynomial $f$.
Note that the triangulated category $D^b{\rm Fuk}^\to(f)$ has a full exceptional collection.

On the other hand, let $f(x,y,z)$ be a weighted homogeneous polynomial. Then one can consider as an analogue of the bounded derived category of coherent sheaves on a smooth proper algebraic variety the following triangulated category. Denote by $S$ the polynomial ring $\CC[x,y,z]$. Let $R_f:= S/(f)$ be the coordinate ring and $L_f$ the maximal grading of $f$. D.~Orlov \cite{Or} considered the triangulated category of a maximally-graded singularity $D_{\rm Sg}^{L_f}(R_f)$ (introduced before by R.-O.~Buchweitz \cite{Bu}) defined as the quotient of the bounded derived category of the category of finitely generated $L_f$-graded $R_f$-modules by the full triangulated subcategory corresponding to finitely generated projective $L_f$-graded $R_f$-modules. It is equivalent to the stable category of $L_f$-graded maximal Cohen-Macaulay modules over $R_f$ and to the stable homotopy category ${\rm HMF}_S^{L_f}(f)$ of $L_f$-graded matrix factorizations of $f$ (see also \cite{Sa3}).

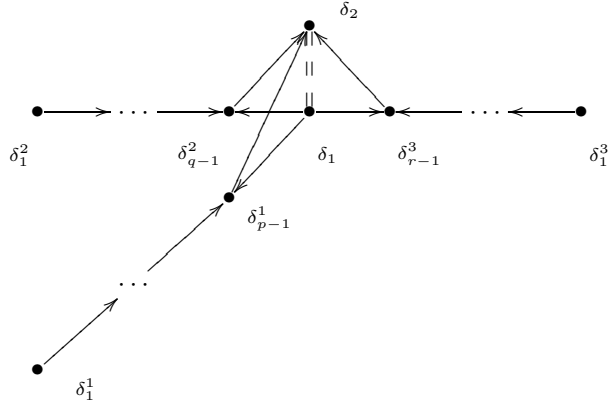
\begin{figure}
$$
\xymatrix{ 
 & & & *{\bullet} \ar@{==}[d] \ar@{<-}[dr]  \ar@{<-}[ldd] \ar@{}^{\delta_{2}}[r]
 & & &  \\
 *{\bullet} \ar@{->}[r] \ar@{}_{\delta^2_1}[d]  & {\cdots} \ar@{->}[r]  & *{\bullet} \ar@{<-}[r] \ar@{->}[ur]   \ar@{}_{\delta^2_{q-1}}[d] & *{\bullet} \ar@{->}[dl] \ar@{->}[r] \ar@{}^{\delta_{1}}[d] & *{\bullet} \ar@{<-}[r]  \ar@{}^{\delta^3_{r-1}}[d]  & {\cdots} \ar@{<-}[r]  &*{\bullet} \ar@{}^{\delta^3_1}[d]   \\
 & &   *{\bullet} \ar@{<-}[dl] \ar@{}_{\delta^1_{p-1}}[r]  & & & & \\
 & {\cdots} \ar@{<-}[dl] & & & & & \\
*{\bullet}  \ar@{}_{\delta^1_1}[r] & & & & & &
  } 
$$
\caption{The quiver $\vec{T}_{p,q,r}$} \label{FigTpqr}
\end{figure}

Moreover, one can consider the quotient stack
${\mathcal C}_f=\left[f^{-1}(0)\backslash\{0\}\left/\widehat{G}_f \right.\right]$ as in Sect.~\ref{sect:Arnold}.
This is a smooth projective line $\PP^1$ with at most three isotropic points of orders $p,q,r$ \cite[Theorem~3]{ET1}. It corresponds to a weighted projective line $\PP^1_{p,q,r}$ with weights $p,q,r$ \cite{GL}. Let $\vec{T}_{p,q,r}$ be the quiver  of Fig.~\ref{FigTpqr} where the double dashed line corresponds to two relations as follows. Let $\beta_1$, $\beta_2$ and $\beta_3$ be the path from $\delta_1$ to $\delta_2$ via $\delta^1_{p-1}$, $\delta^2_{q-1}$ and $\delta^3_{r-1}$ respectively.  Then we consider the relations $\beta_2+\beta_3=0$ and $\beta_1=\beta_3$. They generate an ideal $I$ in the path algebra $\CC \vec{T}_{p,q,r}$ of the quiver. We consider the category ${\rm mod}\mbox{-}\CC \vec{T}_{p,q,r}/I$ of finitely generated right modules over the factor algebra $\CC \vec{T}_{p,q,r}/I$ and its bounded derived category $D^b({\rm mod}\mbox{-}\CC \vec{T}_{p,q,r}/I)$.

Let $D^b{\rm coh}(\PP^1_{p,q,r})$ be the bounded derived category of the category of  coherent sheaves on $\PP^1_{p,q,r}$. W.~Geigle and H.~Lenzing (\cite{GL}, for the special form of the quiver see also \cite[3.9]{LP}) proved the following theorem:
\begin{theorem}[Geigle,Lenzing]
There exists a triangulated equivalence
\[ D^b{\rm coh}(\PP^1_{p,q,r}) \simeq D^b({\rm mod}\mbox{-}\CC \vec{T}_{p,q,r}/I).
\]
\end{theorem}

One has the following $L_f$-graded generalization of Orlov's semi-orthogonal decomposition theorem  \cite[Theorem~2.5]{Or} (see also \cite{t:2}):
\begin{theorem}[Orlov] \label{thm:sod}
\begin{enumerate}
\item If $a_f<0$, one has the semi-orthogonal decomposition
\[
D^b{\rm coh}(\PP^1_{p,q,r}) \simeq \langle D_{\rm Sg}^{L_f}(R_f), {\mathcal A}(0), \ldots , {\mathcal A}(-a_f-1) \rangle,
\]
where ${\mathcal A}(i):= \langle \calO_{\PP^1_{p,q,r}}(\vec{l}) \rangle_{\deg(\vec{l})=i}$.
\item If $a_f=0$, $D^b{\rm coh}(\PP^1_{p,q,r}) \simeq D_{\rm Sg}^{L_f}(R_f)$.
\item If $a_f >0$, one has the semi-orthogonal decomposition
\[ D_{\rm Sg}^{L_f}(R_f) \simeq \langle D^b{\rm coh}(\PP^1_{p,q,r}), {\mathcal K}(0), \ldots , {\mathcal K}(a_f-1) \rangle,
\]
where ${\mathcal K}(i) := \langle (R_f/\mathfrak{m}_f)(\vec{l}) \rangle_{\deg(\vec{l})=i}$ and $\mathfrak{m}_f$ is the maximal ideal in $R_f$.
\end{enumerate}
\end{theorem}

On the other hand, consider a polynomial 
\[ x^p+ y^q + z^r +axyz, \quad \mbox{for some } a \in \CC, a \neq 0. \]
This is called a polynomial of type $T_{p,q,r}$. 
For a triple $(a,b,c)$ of positive integers we define
\[ \Delta(a,b,c) := abc - bc - ac -ab. \]
If $\Delta(p,q,r) >0$ then this polynomial has a cusp singularity at the origin. If $\Delta(p, q, r) =0$ and $a$ is general, then this polynomial has a simple elliptic singularity at the origin. If $\Delta(p, q, r) <0$ then there are other singularities outside the origin and we consider this polynomial as a global polynomial. A distinguished basis of vanishing cycles of such a polynomial (in the case $\Delta(p, q, r) <0$ taking the other singularities into account as well) is given by 
\[(\delta_1, \delta^1_1, \ldots \delta^1_{p-1}, \delta^2_1, \ldots, \delta^2_{q-1}, \delta^3_1, \ldots, \delta^3_{r-1}, \delta_2)\]
with a Coxeter-Dynkin diagram corresponding to the undirected graph $T_{p,q,r}$ underlying the quiver in Fig.~\ref{FigTpqr}.

It is known that the Berglund--H\"{u}bsch duality for some polynomials with nice properties gives the systematic construction 
of mirror pairs of Calabi--Yau manifolds and 
induces the homological mirror symmetry. Therefore, we may expect 
that the homological mirror symmetry can also be categorified to the following:

\begin{sloppypar}
\begin{conjecture}[Takahashi \cite{t:1,t:2}] \label{Conj1}
Let $f(x,y,z)$ be an invertible polynomial.
\begin{enumerate}
\item There should exist a triangulated equivalence
\begin{equation}\label{hms:1}
D^{L_f}_{Sg}(R_f)\simeq  D^b{\rm Fuk}^\to(\widetilde{f}).
\end{equation}
\item There should exist a triangulated equivalence
\begin{equation}\label{hms:2}
D^b{\rm coh}(\PP^1_{p,q,r})\simeq
D^b{\rm Fuk}^\to (T_{p,q,r}).
\end{equation}
\item These triangulated equivalences should be compatible in the following sense: There should exist a diagram
\[ 
\xymatrix{
 D^{L_f}_{Sg}(R_f) \ar[r]^(.4){\sim} \ar@{<->}[d]  & D^b{\rm Fuk}^\to(\widetilde{f}) \ar@{<->}[d]\\
D^b{\rm coh}(\PP^1_{p,q,r}) \ar[r]^(.4){\sim} & D^b{\rm Fuk}^\to (T_{p,q,r})
}
\]
where $D^{L_f}_{Sg}(R_f)$ and $D^b{\rm coh}(\PP^1_{p,q,r})$ are related by Theorem~\ref{thm:sod} and $ D^b{\rm Fuk}^\to(\widetilde{f})$ and $D^b{\rm Fuk}^\to (T_{p,q,r})$ should also be related by semi-orthogonal decomposition.
\end{enumerate}
\end{conjecture}
\end{sloppypar}

A proof of the first part of this conjecture for the simple (ADE) singularities can be derived from a theorem of H.~Kajiura, K.~Saito and A.~Takahashi \cite[Theorem~3.1]{KST1} and results of P.~Seidel \cite[Proposition~3.4]{Se2}, \cite{Se4} (see also \cite{Se3}). Moreover, it was proved by K.~Ueda \cite{Ue} for simple elliptic singularities. The first part of this conjecture can also be stated for invertible polynomials in any number of variables. In this form, it was proved by M.~Futaki and K.Ueda for Brieskorn-Pham singularities \cite{FU1} and for singularities of type D \cite{FU2}. In all these cases, the polynomial $f$ is self-dual, i.e.\ $\widetilde{f}=f$. The second part of this conjecture was proved for the case $r=1$ by P.~Seidel \cite{Se2}, D.~van Straten \cite{vSt} and D.~Auroux, L.~Katzarkov and D.~Orlov \cite{AKO}, for the case $r=2$ by A.~Takahashi \cite{t:3} and in general by A.~Keating \cite{Keating}.

Now consider an invertible polynomial $f(x,y,z)$  and a finite group $G$ of diagonal symmetries of $f$. We assume that $G$ contains the group $G_0$ generated by the exponential grading operator $g_0$. The orbifold curve ${\mathcal C}_{(f,G)}:=\left[f^{-1}(0)\backslash\{0\}\left/\widehat{G} \right.\right]$ is mirror dual to the following data: A function $F : U \to \CC$, defined on a suitably chosen submanifold $U$ of $\CC^3$, given by $F(x,y,z)=x^{\gamma_1'} + y^{\gamma_2'} + z^{\gamma_3'}-xyz$. The group $\widetilde{G}$ leaves $F$ invariant and we can consider a crepant resolution $Y \to U/\widetilde{G}$ given by the $\widetilde{G}$-Hilbert scheme and the 
proper transform $\widehat{X} \subset Y$ of $X=F^{-1}(0)/\widetilde{G} \subset U/\widetilde{G}$ (cf.\ \cite{Se5}).

Let ${\rm HMF}_S^{\widehat{G}}(f)$ be the stable homotopy category of $\widehat{G}$-graded matrix factorizations of $f$. Let $D^b{\rm Coh}{\mathcal C}_{(f,G)}$ be the derived category of the category of coherent sheaves on ${\mathcal C}_{(f,G)}$.

We arrive at the following generalization of Conjecture~\ref{Conj1} (cf.\  \cite{ET2}):
\begin{conjecture}[E., Takahashi]  \label{Conj2}
There should exist triangulated equivalences 
\[ 
\xymatrix{
 {\rm HMF}_S^{\widehat{G}}(f) \ar[r]^(.4){\sim} \ar@{<->}[d]  & D^b{\rm Fuk}^\to(\widetilde{f})//\widetilde{G} \ar@{<->}[d]\\
D^b{\rm Coh}{\mathcal C}_{(f,G)} \ar[r]^(.4){\sim} & D^b{\rm Fuk}^\to(F)//\widetilde{G}
}
\]
where the two lines are related by semi-orthogonal decompositions,
$F(x,y,z)=x^{\gamma_1'} + y^{\gamma_2'} + z^{\gamma_3'}-xyz$ is right equivalent  to $\widetilde{f}(x,y,z)-xyz$, and $-//\widetilde{G}$ means the smallest triangulated category containing 
the orbit category $-/\widetilde{G}$ $($cf.\ \cite{Asashiba:2008, Cibils_Marcos:2006} for orbit categories; see also \cite{Keller:2005}$)$.
\end{conjecture}

\section{Strange duality} \label{sect:SD}
We now give some evidence for the conjectures stated in the last section.

Let $f(x_1,\dots ,x_n)$ be an invertible polynomial and $G \subset G_f$ a subgroup of the maximal group of symmetries. We shall investigate the correspondence
\[ (f,G) \longleftrightarrow (\widetilde{f}, \widetilde{G}). 
\]
First let $n=3$, $f(x,y,z)$ be a non-degenerate invertible polynomial such that $\widetilde{f}(x,y,z)$ is non-degenerate as well and let $G=G_f$. Then the correspondence
\[ (f,G_f) \longleftrightarrow (\widetilde{f}, \{ e \})
\]
was considered in \cite{ET1}. We defined Dolgachev numbers for a pair $(f,G_f)$ and Gabrielov numbers for a pair $(f, \{ e \})$ as follows.

The quotient stack
\[
{\mathcal C}_{(f,G_f)}:=\left[f^{-1}(0)\backslash\{0\}\left/\widehat{G}_f\right.\right]
\]
is a smooth projective line $\PP^1$ with at most three isotropic points of orders $\alpha_1, \alpha_2, \alpha_3$ (see Sect.~\ref{sect:HMS}). 
\begin{definition}
The numbers $(\alpha_1, \alpha_2, \alpha_3)$ are called the {\em Dolgachev numbers} of the pair $(f,G_f)$ and the tuple is denoted by $A_{(f, G_f)}$. 
\end{definition}

On the other hand, consider
the deformation $F(x,y,z):=f(x,y,z)-xyz$ of $f$. By \cite[Theorem~10]{ET1}, if $\Delta(\gamma_1, \gamma_2, \gamma_3) >0$ there exists a holomorphic coordinate change so that this polynomial becomes a polynomial of type $T_{\gamma_1, \gamma_2, \gamma_3}$ (for the definition see Sect.~\ref{sect:HMS}). In the cases $\Delta(\gamma_1, \gamma_2, \gamma_3) =0$ and $\Delta(\gamma_1, \gamma_2, \gamma_3) <0$ there is also a relation to a polynomial of type $T_{\gamma_1, \gamma_2,\gamma_3}$, see \cite[Theorem~10]{ET1}.
\begin{definition}
The numbers $(\gamma_1, \gamma_2, \gamma_3)$ are called the {\em Gabrielov numbers} of the pair $(f, \{ e \})$ and the tuple is denoted by $\Gamma_{(f, \{ e \})}$.
\end{definition}

By \cite[Theorem~13]{ET1} we have the following theorem:

\begin{theorem}[E., Takahashi]  \label{thm:ET1}
Let $f(x,y,z)$ be a non-degenerate invertible polynomial such that $\widetilde{f}(x,y,z)$ is non-degenerate as well. Then we have 
\[ A_{(f, G_f)} = \Gamma_{(\widetilde{f}, \{ e \})}, \qquad A_{(\widetilde{f}, G_{\widetilde{f}})}= \Gamma_{(f, \{ e \})}.
\]
Namely, the Dolgachev numbers  $A_{(f, G_f)}$ for the pair $(f,G_f)$ coincide with the Gabrielov numbers $\Gamma_{(\widetilde{f}, \{ e \})}$ for the pair $(\widetilde{f}, \{ e \})$ and the Dolgachev numbers $A_{(\widetilde{f}, G_{\widetilde{f}})}$ for the pair $(\widetilde{f}, G_{\widetilde{f}})$ coincide with the Gabrielov numbers $\Gamma_{(f, \{ e \})}$ for the pair $(f, \{ e \})$.
\end{theorem}

The 14 exceptional unimodal singularities can be given by non-degenerate invertible polynomials $f(x,y,z)$ with $G_f=G_0$. These are the polynomials indicated in Table~\ref{TabArnold}. The Dolgachev and Gabrielov numbers coincide with the corresponding numbers indicated in this table. Therefore we obtain Arnold's strange duality as a special case of this theorem. We come back to this duality in Sect.~\ref{sect:Ex}.

More generally, let $f(x,y,z)$ be a non-degenerate invertible polynomial such that $\widetilde{f}(x,y,z)$ is non-degenerate as well, but now consider a subgroup $G$ with $G_0 \subset G \subset G_f$. Then $\{ e \} \subset \widetilde{G} \subset {\rm SL}_n(\ZZ) \cap G_{\widetilde{f}}$. In \cite{ET2} we defined Dolgachev numbers for the pair $(f,G)$ with $G_0 \subset G$ and Gabrielov numbers for a pair $(f,G)$ with $G \subset {\rm SL}_n(\ZZ)$ as follows.

The quotient stack
\[
{\mathcal C}_{(f,G)}:=\left[f^{-1}(0)\backslash\{0\}\left/\widehat{G}\right.\right]
\]
can be regarded as a smooth projective curve of genus $g_{(f,G)}$ with a finite number of isotropic points. 
\begin{definition}
The orders $\alpha_1, \ldots, \alpha_r$ of the isotropy groups of these points are called the {\em Dolgachev numbers} of the pair $(f,G)$ and denoted by $A_{(f, G)}$.
\end{definition}

By \cite[Theorem~4.6]{ET2}, the Dolgachev numbers of the pair $(f,G)$ can be computed from the Dolgachev numbers of the pair $(f,G_f)$ as follows. Let $A_{(f,G_f)} = (\alpha_1', \alpha_2', \alpha_3')$ be the Dolgachev numbers of the pair $(f,G_f)$.
For positive integers $u$ and $v$, by $u \ast v$ we denote $v$ copies of the integer $u$. 

\begin{theorem}[E., Takahashi]  \label{thm:Dol}
Let $H_i \subset G_f$ be the minimal subgroup containing $G$ and the isotropy group of the point $p_i$, $i=1,2,3$. Then we have the following formula for the Dolgachev numbers $\alpha_1, \ldots, \alpha_r$ of the pair $(f,G)$:
\[ (\alpha_1, \ldots, \alpha_r) = \left( \frac{\alpha'_i}{|H_i/G|} \ast |G_f/H_i|,\ i=1,2,3 \right) ,\]
where we omit numbers which are equal to one on the right-hand side.
\end{theorem}

We define the {\em stringy Euler number} of the orbifold curve ${\mathcal C}_{(f,G)}$ by
\[ e_{\rm st}({\mathcal C}_{(f,G)}) := 2-2g_{(f,G)} +\sum_{i=1}^r (\alpha_i-1). \]

Now consider a pair $(f,G)$ with $G \subset {\rm SL}_3(\ZZ)$. 
\begin{definition}
Let $\Gamma_{(f, \{ e \})}= (\gamma'_1,\gamma'_2, \gamma'_3)$ be the Gabrielov numbers of the pair $(f, \{ e \})$ and let $K_i \subset G$ be the maximal subgroup of $G$ fixing the coordinate $x_i$, $i=1,2,3$. Then the {\em Gabrielov numbers} of the pair $(f,G)$ are the numbers $\gamma_1, \dots, \gamma_s$ defined by 
\[ (\gamma_1, \ldots, \gamma_s) = \left( \frac{\gamma'_i}{|G/K_i|} \ast |K_i|,\ i=1,2,3 \right) ,\]
where we omit numbers which are equal to one on the right-hand side. We denote this tuple of numbers by $\Gamma_{(f,G)}$. 
\end{definition}
In \cite{ET4}, we gave a geometric definition of these numbers as lengths of arms of a certain Coxeter-Dynkin diagram: Let $U$ be a suitably chosen submanifold of $\CC^3$. We consider a crepant resolution $Y \to U/G$ and the preimage $Z$ of the image of the Milnor fibre of the cusp singularity $T_{\gamma'_1, \gamma'_2, \gamma'_3}$ under the natural projection $U \to U/G$. Using the McKay correspondence, we constructed a basis of the relative homology group $H_3(Y,Z;\QQ)$ with a Coxeter-Dynkin diagram where one can read off the Gabrielov numbers.

Let $G$ be a finite group acting linearly on $\CC^n$. For an element $g\in G$, its {\em age} \cite{Ito-Reid}
is defined by ${\rm age\,}(g):=\sum_{i=1}^n\alpha_i$,
where in a certain basis in $\CC^n$ one has
$g={\rm diag\,}({\bf e}[\alpha_1], \ldots, {\bf e}[\alpha_n])$
with $0\le\alpha_i < 1$. Now let $G \subset {\rm SL}_n(\ZZ)$. Then the age of an element $g \in G$ is an integer. Define
\[ j_{G}:= \{ g \in G \, | \, {\rm age}(g)=1,\ g \mbox{ only fixes the origin} \}.
\]

Let  $F$ be a polynomial of type $T_{\gamma'_1, \gamma'_2, \gamma'_3}$ with 
the Gabrielov numbers $(\gamma'_1,\gamma'_2,\gamma'_3)$ for the pair $(f,G)$,
Define the {\em $G$-equivariant Milnor number} of $F$ by
\[
\mu_{(F,G)} :=  2- 2j_G + \sum_{i=1}^s (\gamma'_i -1). 
\]

By \cite[Theorem~7.1]{ET2} we have the following result:

\begin{theorem}[E., Takahashi]  \label{thm:ET2}
Let $f(x,y,z)$ be a non-degenerate invertible polynomial such that $\widetilde{f}(x,y,z)$ is non-degenerate as well and let $G_0 \subset G \subset G_f$. Then we have
\[ g_{(f,G)} = j_{\widetilde{G}}, \quad A_{(f,G)}=\Gamma_{(\widetilde{f},\widetilde{G})}, \quad 
 e_{\rm st}({\mathcal C}_{(f,G)}) = \mu_{(F,\widetilde{G})},
 \]
where $F$ is a polynomial of type $T_{\gamma'_1, \gamma'_2, \gamma'_3}$ with 
the Gabrielov numbers $(\gamma'_1,\gamma'_2,\gamma'_3)$ for the pair $(\widetilde{f},\{ e\})$.
\end{theorem}

\section{Orbifold E-functions} \label{sect:E}
We now show that the mirror symmetry for Berglund-H\"ubsch dual pairs also holds on the level of suitably defined Hodge numbers. Therefore we discuss the notion of an orbifold E-function for a polynomial with an isolated singularity at the origin.

Let  $f(x_1,\dots, x_n)$ be a polynomial with $f(0)=0$ and with an isolated singularity at 0. 
We regard the polynomial $f$ as a holomorphic map $f:V\to\CC$ where 
$V$ is a suitably chosen neighbourhood of $0\in\CC^n$ so that the fibration 
$f$ has good technical properties.
Consider the Milnor fibre $V_f:=\{x\in V~|~f(x)=1\}$ of the fibration $f:X\to\CC$. 
J.~H.~M.~Steenbrink \cite{st:1} constructed a canonical mixed Hodge structure  on the vanishing cohomology $H^{n-1}(V_f,\CC)$ 
with an automorphism $c$ 
given by the Milnor monodromy. 
We can naturally associate a bi-graded vector space to a mixed Hodge structure 
with an automorphism.
Consider the Jordan decomposition $c=c_{\rm ss}\cdot c_{\rm unip}$ of $c$ where $c_{\rm ss}$ and $c_{\rm unip}$ denote 
the semi-simple part and unipotent part respectively.
For $\lambda\in \CC$, let  
\begin{equation}
H^{n-1}(V_f,\CC)_\lambda:={\rm Ker}(c_{\rm ss}-\lambda\cdot {\rm id}:H^{n-1}(V_f,\CC)\longrightarrow 
H^{n-1}(V_f,\CC)).
\end{equation}
Denote by $F^\bullet$ the Hodge filtration of the mixed Hodge structure.

\begin{definition}
Define the $\QQ\times \QQ$-graded vector space $\H_f:=\displaystyle\bigoplus_{p,q\in\QQ}\H^{p,q}_f$ as
\begin{enumerate}
\item If $p+q\ne n$, then  $\H^{p,q}_f:=0$.
\item If $p+q=n$ and $p\in\ZZ$, then  
\[
\H^{p,q}_f:={\rm Gr}^{p}_{F^\bullet}H^{n-1}(V_f,\CC)_1.
\]
\item If $p+q=n$ and $p\notin\ZZ$, then  
\[
\H^{p,q}_f:={\rm Gr}^{[p]}_{F^\bullet}H^{n-1}(V_f,\CC)_{e^{2\pi\sqrt{-1} p}},
\]
where $[p]$ is the largest integer less than $p$.
\end{enumerate}
\end{definition}

Let $G$ be a subgroup of the maximal group $G_f$ of diagonal symmetries of $f$. For $g \in G$, we denote by ${\rm Fix}\, g :=\{x\in\CC^n~|~g\cdot x=x \}$ the fixed locus of $g$, 
by $n_g: = \dim {\rm Fix}\, g$ its dimension and by $f^g:=f|_{{\rm Fix}\, g}$ the restriction of $f$ to the fixed locus of $g$. 
Note that the function $f^g$ has an isolated singularity at the origin \cite[Proposition~5]{ET3}.

We shall use the fact that $\H_{f^g}$ admits a natural $G$-action 
by restricting the $G$-action on $\CC^n$ to ${\rm Fix}\, g$ (which is well-defined since $G$ acts diagonally on $\CC^n$).

To the pair $(f,G)$ we can associate the following $\QQ\times \QQ$-graded super vector space:
\begin{definition}
Define the $\QQ\times \QQ$-graded super vector space $\H_{f,G}:=\H_{f,G,\bar{0}}\oplus \H_{f,G,\bar{1}}$ as 
\begin{equation}
\H_{f,G,\overline{0}}:=\bigoplus_{\substack{g\in G;\\ n_g\equiv 0\ (\text{\rm mod } 2)}}(\H_{f^g})^G(-{\rm age}(g),-{\rm age}(g)),
\end{equation}
\begin{equation}
\H_{f,G,\overline{1}}:=\bigoplus_{\substack{g\in G;\\ n_g\equiv 1\ (\text{\rm mod } 2)}}(\H_{f^g})^G(-{\rm age}(g),-{\rm age}(g)),
\end{equation}
where $(\H_{f^g})^G$ denotes the $G$-invariant subspace of $\H_{f^g}$.
\end{definition}

\begin{definition}[\cite{EGT}]
The {\em E-function} for the pair $(f,G)$ is 
\begin{equation}
E(f,G)(t,\bar{t})=\displaystyle\sum_{p,q\in\QQ}\left({\rm dim}_\CC (\H_{f,G,\bar{0}})^{p,q}-{\rm dim}_\CC (\H_{f,G,\bar{1}})^{p,q}\right)
\cdot t^{p-\frac{n}{2}}{\bar{t}}^{q-\frac{n}{2}}.
\end{equation}
\end{definition}
In general, we may have both $(\H_{f,G,\overline{0}})^{p,q}\ne 0$ and $(\H_{f,G,\overline{1}})^{p,q}\ne 0$ for some $p,q\in\QQ$ (see \cite{EGT}).
However we have the following proposition (see \cite[Proposition~3]{EGT}):
\begin{proposition}
Let $f(x_1, \ldots , x_n)$ be a non-degenerate invertible polynomial and $G$ a subgroup of $G_f$.
Assume $G\subset {\rm SL}(n;\CC)$ or $G\supset G_0$.
If $(\H_{f,G,\overline{i}})^{p,q}\ne 0$, then $(\H_{f,G,\overline{i+1}})^{p,q}= 0$ for all $p,q\in\QQ$ and $\overline{i}\in\ZZ/2\ZZ$.
\end{proposition}

\begin{definition}
Let $f(x_1, \ldots , x_n)$ be a non-degenerate invertible polynomial and $G$ a subgroup of $G_f$.
Assume $G\subset {\rm SL}(n;\CC)$ or $G\supset G_0$.
The {\em Hodge numbers} for the pair $(f,G)$ are
\[
h^{p,q}(f,G):=\dim_\CC (\H_{f,G,\overline{0}})^{p,q}+\dim_\CC (\H_{f,G,\overline{1}})^{p,q},\quad p,q\in\QQ.
\]
\end{definition}

\begin{proposition}\label{prop4}
Let $f(x_1, \ldots , x_n)$ be a non-degenerate invertible polynomial and $G$ a subgroup of $G_f$.
The E-function is given by 
\[
E(f,G)(t,\bar{t}):=
\begin{cases}
\displaystyle\sum_{p,q\in\QQ}(-1)^{p+q}h^{p,q}(f,G) \cdot 
t^{p-\frac{n}{2}}\bar{t}^{q-\frac{n}{2}},\ \text{if }G \subset {\rm SL}_n(\CC),\\
\displaystyle\sum_{p,q\in\QQ}(-1)^{-p+q}h^{p,q}(f,G) \cdot 
t^{p-\frac{n}{2}}\bar{t}^{q-\frac{n}{2}}, \ \text{if } G_0 \subset G.
\end{cases}
\]
\end{proposition}

Therefore, in the case that $f$ is a non-degenerate invertible polynomial and $G \subset {\rm SL}_n(\CC)$, the definition of the E-function for the pair $(f,G)$ agrees with \cite[Definition~5.7]{ET2}:

\begin{definition} Let $f(x_1, \ldots , x_n)$ be a polynomial with an isolated singularity at the origin invariant under a group $G \subset {\rm SL}_n(\CC)$. The {\em E-function} of the pair $(f,G)$ is defined by
\[
E(f,G)(t,\bar{t}):=\displaystyle\sum_{p,q\in\QQ}(-1)^{p+q}h^{p,q}(f,G) \cdot 
t^{p-\frac{n}{2}}\bar{t}^{q-\frac{n}{2}}.
\]
\end{definition}
The E-function of the pair $(f,G)$ is the generating function of the exponents of the pair $(f,G)$. An {\em exponent} of the pair $(f,G)$ is a rational number $q$ with $h^{p,q}(f,G)\ne 0$. 
The {\em set of exponents} of the pair $(f,G)$ is the multi-set of exponents 
\[
\left\{q*h^{p,q}(f,G)~|~p,q\in\QQ,\ h^{p,q}(f,G)\ne 0 \right\},
\]
where by $u*v$ we denote $v$ copies of the rational number $u$.  
It is almost clear that the mean of the set of exponents of $(f,G)$ is $n/2$, namely, we have 
\[
\sum_{p,q\in\QQ}(-1)^{-p+q} \left(q-\frac{n}{2} \right)h^{p,q}(f,G)=0.
\]
It is natural to ask what is the {\em variance of the set of exponents} of $(f,G)$ defined by 
\[ {\rm Var}_{(f,G)} := \sum_{p,q\in\QQ}(-1)^{-p+q} \left(q- \frac{n}{2} \right)^2 h^{p,q}(f,G). \]
In \cite{ET3} we have proved:
\begin{theorem}[E., Takahashi]  \label{thm:Var}
Let $f(x_1, \ldots, x_n)$ be a non-degenerate weighted homogeneous polynomial invariant under a group $G \subset {\rm SL}_n(\CC)$.
Then one has
\[
{\rm Var}_{(f,G)} = \frac{1}{12} \hat{c} \cdot\chi(f,G),
\]
where $\hat{c} := n - 2\sum_{i=1}^n q_i$ and $\chi(f,G):=E(f,G)(1,1)$.
\end{theorem}

We also define (see \cite[Definition~5.3, Definition~5.10]{ET2}):
\begin{definition} Let $f(x_1, \ldots , x_n)$ be a polynomial with an isolated singularity at the origin invariant under a group $G \subset {\rm SL}_n(\CC)$.
The {\em characteristic polynomial} of the pair $(f,G)$ is 
\[
\phi(f,G)(t):= \prod_{q\in\QQ} (t- {\bf e}[q])^{h^{p,q}(f,G)}.
\]
\end{definition}

We have (see \cite[Theorem~5.12]{ET2}):
\begin{theorem}[E., Takahashi]
Let $F(x_1,x_2,x_3)= x_1^{\gamma'_1} + x_2^{\gamma'_2} + x_3^{\gamma'_3} - x_1x_2x_3$ and $G$ be a subgroup of $G_F$. Then
\[ \phi(F,G)(t) = (t-1)^{2-2j_G} \prod_{i=1}^s \frac{t^{\gamma_i}-1}{t-1}, \]
where $\gamma_1, \ldots, \gamma_s$ are the Gabrielov numbers defined in Sect.~\ref{sect:SD}.
\end{theorem}

The characteristic polynomial $\phi(f,G)(t)$ agrees with the reduced orbifold zeta function $\overline{\zeta}_{f,G}^{\rm orb}(t)$ defined in \cite{EGPEMS}. Its degree is the reduced orbifold Euler characteristic $\overline{\chi}(V_f,G)$  of $(V_f,G)$ (see below). 

From physical reasons \cite{BHe}, one expects that there is the following relation between the orbifold E-functions of dual pairs. This was proved in \cite{EGT}.
\begin{theorem}[E., Gusein-Zade, Takahashi]  \label{main}
Let $f(x_1, \ldots , x_n)$ be a non-degenera\-te invertible polynomial and $G$ a subgroup of $G_f$. Then 
\[
E(f,G)(t, \bar{t}) = (-1)^n  E(\widetilde{f},\widetilde{G})(t^{-1},\bar{t}). 
\]
\end{theorem}

Using Proposition~\ref{prop4}, we can derive from this theorem the mirror symmetry of dual pairs on the level of Hodge numbers:
\begin{corollary}
Let $f(x_1, \ldots , x_n)$ be a non-degenerate invertible polynomial and $G$ a subgroup of $G_f\cap {\rm SL}(n;\CC)$. Then 
for all $p,q\in\QQ$, one has
\[
h^{p,q}(f,G)=h^{n-p,q}(\widetilde{f},\widetilde{G}).
\]
\end{corollary}

As another corollary, we get the main result of \cite{EGPEMS}:
\begin{corollary} \label{cor:zeta}
 One has
 \[
  \overline{\zeta}^{{\rm orb}}_{\widetilde{f},\widetilde{G}}(t)=\left(
  \overline{\zeta}^{{\rm orb}}_{f,G}(t)  
  \right)^{(-1)^n}\,.
 \]
 \end{corollary}

From this we derive the main result of \cite{EGMMJ2}:
\begin{corollary} \label{cor:Euler}
One has 
$$ 
\overline{\chi}(V_{\widetilde{f}}, \widetilde{G})=(-1)^n\overline{\chi}(V_f,G)\,. 
$$ 
\end{corollary}
Note that the latter two results were even proven without the assumption of non-degeneracy.

We also obtain as a corollary from Theorem~\ref{thm:Var} and Theorem~\ref{main}:
\begin{corollary}
Let $f(x_1, \ldots , x_n)$ be a non-degenerate invertible polynomial and $G$ a subgroup of $G_f$ containing $G_0$. Then one has
\[
{\rm Var}_{(f,G)} = \frac{1}{12} \hat{c} \cdot\chi(f,G),
\]
where $\hat{c} := n - 2\sum_{i=1}^n q_i$ and $\chi(f,G):=E(f,G)(1,1)$.
\end{corollary}

\section{Saito duality} \label{sect:Saito}
In this section we consider a generalization of Saito's duality (Theorem~\ref{thm:Saito}) to the Berglund-H\"ubsch duality.

For this we recall the notion of the Burnside ring of a finite group (see \cite{Knutson}). Let $G$ be a finite group. A $G$-set is a set with an action of the group $G$. A $G$-set is {\em irreducible} if the action
of $G$ on it is transitive. Isomorphism classes of irreducible $G$-sets are in one-to-one correspondence with
conjugacy classes of subgroups of $G$: to the conjugacy class containing a subgroup $H\subset G$ one associates the isomorphism class $[G/H]$
of the $G$-set $G/H$. 
\begin{definition} The {\em Burnside ring} $B(G)$ of $G$ is the Grothendieck ring of finite $G$-sets, i.e.\ it is the (abelian) group generated by the isomorphism classes of finite $G$-sets
modulo the relation $[A\amalg B]=[A]+[B]$ for finite $G$-sets $A$ and $B$. The multiplication
in $B(G)$ is defined by the cartesian product.
\end{definition}
As an abelian group, $B(G)$
is freely generated by the isomorphism classes of irreducible $G$-sets. The element $1$ in
the ring $B(G)$ is represented by the $G$-set consisting of one point (with the trivial $G$-action).

Let $V$ be a ``good'' topological space, say, a union of cells in a finite $CW$-complex or a quasi-projective complex or real analytic variety. Let $G$ be a finite group acting on $V$.  For $x \in V$, denote by $G_x$ the isotropy subgroup of $x$. For a subgroup $H \subset G$ let $V^H$ be the set of all fixed points of $H$. Denote by
$V^{(H)}$ the set of points of $V$ with isotropy group $H$. Finally, let 
\[ 
V^{([H])}:= \bigcup_{K \in [H]} V^{(K)}.
\]
The equivariant Euler characteristic was defined in \cite{Verdier, TtD}.
\begin{definition}
The {\em equivariant Euler characteristic} of the $G$-space $V$ is defined by
$$
\chi^G(V):=\sum\limits_{[H]\in{\rm Conjsub\,}G}\chi(V^{([H])}/G)[G/H] \in B(G),
$$
where $\Conjsub\, G$ denotes the set of conjugacy classes of subgroups of $G$ and $\chi(Z)$ denotes the usual Euler characteristic of the topological space $Z$.
\end{definition}

There is also the notion of an orbifold Euler characteristic (\cite{DHVW1,DHVW2}, see also \cite{HH} and the references therein):
\begin{definition}
 The {\em orbifold Euler characteristic} of the $G$-space $V$ is defined by 
$$ 
\chi^{\rm orb}(V,G)=\frac{1}{\vert G\vert}\sum_{(g,h):gh=hg} \chi(X^{\langle g,h\rangle}) 
$$ 
where $\langle g,h\rangle$ is the subgroup of $G$ generated by $g$ and $h$.
\end{definition}

There is a map $r_G^{\rm orb}: B(G)\to \ZZ$ defined by sending a class
$[G/H]$  to the number $\chi^{\rm orb}([G/H],G)$. If $G$ is abelian then $\chi^{\rm orb}([G/H],G)=|H|$. We have
\[
r_G^{\rm orb}(\chi^G(V))=\chi^{\rm orb}(V,G).
\]

\begin{definition}
\begin{itemize}
\item The {\em reduced equivariant Euler characteristic} of the $G$-space $V$ is
\[
\overline{\chi}^G(V) := \chi^G(V) -1.
\]
\item The {\em reduced orbifold Euler characteristic} of the $G$-space $V$ is
\[
\overline{\chi}^{\rm orb}(V,G)= \chi(V,G)- \vert G\vert. 
\]
\end{itemize}
\end{definition}
We have
$$
r_G^{\rm orb}(\overline{\chi}^G(V))=\overline{\chi}^{\rm orb}(V,G).
$$

For a group $G$ let $G^\ast:= {\rm Hom}(G, \CC^\ast)$ be its group of characters.
\begin{definition}
Let $G$ be abelian.
The {\em equivariant Saito duality} between $B(G)$ and $B(G^\ast)$ is the map
\[
\begin{array}{cccc}
D_G: & B(G) &  \to &  B(G^\ast) \\
              & [G/H]  & \mapsto & [G^\ast/\widetilde{H}]        
\end{array}
\]
\end{definition}

In \cite{EGBLMS} it was proved:
\begin{theorem}[E., Gusein-Zade] \label{thm:Burn}
Let $f(x_1, \ldots ,x_n)$ be an invertible polynomial. Then one has
\[ \overline{\chi}^{G_{\widetilde{f}}}(V_{\widetilde{f}}) = (-1)^n D_{G_f}(\overline{\chi}^{G_f}(V_f))
\]
\end{theorem}   

In the special case when the groups of diagonal symmetries of the dual polynomials
are cyclic and are generated by the monodromy transformations, this yields the original Saito duality (Theorem~\ref{thm:Saito}). Moreover, it is shown in \cite{EGBLMS} that the relation between ``geometric roots''  of the zeta functions for some Berglund-H\"ubsch dual invertible polynomials described in \cite{EGMMJ1} is a special case of Theorem~\ref{thm:Burn}.         
One can also derive Corollary~\ref{cor:Euler} from this theorem.

In order to derive Corollary~\ref{cor:zeta} from a similar result, we considered in \cite{EG1506} an enhancement of the Burnside ring:
\begin{definition} \label{def1}
A {\em finite enhanced} $G$-{\em set} is a triple $(X, h, \alpha)$, where:
\begin{enumerate}
 \item[1)] $X$ is a finite $G$-set;
 \item[2)] $h$ is a one-to-one $G$-equivariant map $X\to X$;
 \item[3)] $\alpha$ associates to each point $x\in X$ a one-dimensional (complex) representation
 $\alpha_x$ of the isotropy subgroup $G_x=\{a\in G: ax=a\}$ of the point $x$ so that:
 \begin{itemize}
 \item[(a)] for $a\in G$ one has $\alpha_{ax}(b)=\alpha_x(a^{-1}ba)$, where $b\in G_{ax}=aG_xa^{-1}$;
 \item[(b)] $\alpha_{h(x)}(b)=\alpha_x(b)$.
 \end{itemize}
\end{enumerate}
\end{definition}
\begin{definition} 
The {\em enhanced Burnside ring} $\widehat{B}(G)$ is the Grothendieck group of finite enhanced $G$-sets.
\end{definition}

Let $V$ be a complex manifold with a complex analytic action of a finite group $G$. 
For a point $x \in V$ we consider the one-dimensional representation $\alpha_x : G_x \to \CC^\ast$ defined by $\alpha_x(g) = {\bf e}[\mbox{age}(g)]$.
Let $\varphi: V \to V$ be a $G$-equivariant map with $\alpha_{\varphi(x)}=\alpha_x$ for all $x \in V$. In \cite{EG1506} we defined an
{\em enhanced Euler characteristic}  $\widehat{\chi}^G(V,\varphi) \in \widehat{B}(G)$ and a {\em reduced enhanced Euler characteristic} $\overline{\widehat{\chi}}^G(V,\varphi)$ of the pair $(V,\varphi)$.

Now let $G$ again be abelian.
Let $\widehat{B}_1(G)$ be the subgroup of $\widehat{B}(G)$ generated by the isomorphism classes of finite enhanced sets $(X,h,\alpha)$ such that $h(x)\in G x$ for all $x\in X$.
As an abelian group it is freely generated by the classes $[G/H,h, \alpha]$ where 
\begin{enumerate}
\item $h: G/H \to G/H$ can be identified with an element $h \in G/H$,
\item $\alpha$ is a 1-dimensional representation of $H$. 
\end{enumerate}
The factor group $G/H$ is canonically isomorphic to $\widetilde{H}^\ast=\text{Hom}(\widetilde{H},\CC^*)$
and the group of characters $H^*=\text{Hom}(H,\CC^*)$ is canonically isomorphic to $G^*/\widetilde{H}$. In this way, the element $h \in G/H$ defines a 1-dimensional representation $\widetilde{h}: \widetilde{H} \to \CC^\ast$ and the representation $\alpha: H \to \CC^\ast$ defines an element $\widetilde{\alpha} \in G^\ast/\widetilde{H}$.

\begin{definition}
The {\em enhanced equivariant Saito duality} between $\widehat{B}_1(G)$ and $\widehat{B}_1(G^\ast)$ is the map
\[
\begin{array}{cccc}
\widehat{D}_G: & \widehat{B}_1(G) &  \to &  \widehat{B}_1(G^\ast) \\
              & [G/H,h,\alpha]  & \mapsto & [G^\ast/\widetilde{H},\widetilde{\alpha},\widetilde{h}]        
\end{array}
\]
\end{definition}

In \cite{EG1506} we proved:
\begin{theorem}[E., Gusein-Zade]
Let $f(x_1, \ldots ,x_n)$ be an invertible polynomial and let $h_f: V_f \to V_f$ and $h_{\widetilde{f}} : V_{\widetilde{f}} \to V_{\widetilde{f}}$ be the monodromy transformations of $f$ and $\widetilde{f}$ respectively. Then one has
 $$
 \overline{\widehat{\chi}}^{G_{\widetilde{f}}}(V_{\widetilde{f}},h_{\widetilde{f}})=
 (-1)^n \widehat{D}_{G_f}(\overline{\widehat{\chi}}^{G_f}(V_f,h_f))
 $$
\end{theorem}

It is shown in \cite{EG1506} that one can derive Corollary~\ref{cor:zeta} from this theorem.

\section{Examples} \label{sect:Ex}
We first consider Arnold's classification of singularities \cite{Ar}.

\begin{figure}[h]
$$
\xymatrix{ 
  *{\bullet} \ar@{-}[r] \ar@{}_{E^2_1}[d]  & {\cdots} \ar@{-}[r]  & *{\bullet} \ar@{-}[r]  \ar@{}_{E^2_{\alpha_2-1}}[d] & *{\bullet} \ar@{-}[dl] \ar@{-}[r] \ar@{}^{E}[d] & *{\bullet} \ar@{-}[r]  \ar@{}^{E^3_{\alpha_3-1}}[d]  & {\cdots} \ar@{-}[r]  &*{\bullet} \ar@{}^{E^3_1}[d]   \\
 & &   *{\bullet} \ar@{-}[dl] \ar@{}_{E^1_{\alpha_1-1}}[r]  & & & & \\
 & {\cdots} \ar@{-}[dl] & & & & & \\
*{\bullet}  \ar@{}_{E^1_1}[r] & & & & & &
  } 
$$
\caption{The dual graph of $\mathcal E$} \label{FigExc}
\end{figure}
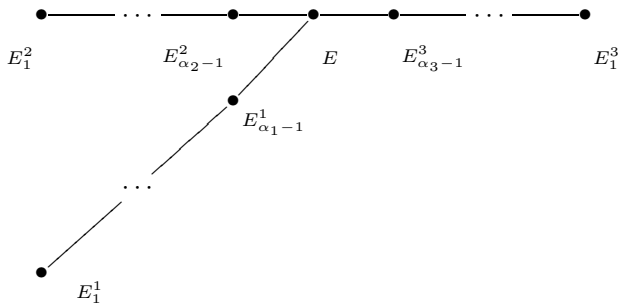

We first have the simple singularities which are also called the ADE singularities. They are given by invertible polynomials $f$ with $a_f<0$. These polynomials together with the corresponding Dolgachev and Gabrielov numbers are given in \cite[Table~5]{ET1}. 
\begin{table}[h]
\begin{center}
\begin{tabular}{cccc}
\hline
{\rm Name} & $f(x,y,z)$ & $\alpha_1, \alpha_2, \alpha_3$ & $\gamma_1, \gamma_2, \gamma_3$ \\
\hline
$A_k$ & $xy+y^kz+zx$, $k \geq 1$ & $k,1,1$ & $1,1,k$ \\
$D_{2k+1}$    & $x^2+xy^k+yz^2$, $k \geq 2$ & $2,2,2k-1$ & $2,2,2k-1$  \\
$D_{2k+2}$ & $x^2+y^2z+yz^{k+1}$, $k \geq 1$ & $2,2,2k$ &  $2,2,2k$ \\
$E_6$  & $x^3+y^2+yz^2$ & $3,2,3$ & $3,2,3$   \\
$E_7$  & $x^2+y^3+yz^3$ & $2,3,4$ & $2,3,4$  \\
$E_8$  & $x^2+y^3+z^5$ & $2,3,5$ & $2,3,5$  \\ 
\hline
\end{tabular}
\end{center}
\caption{The ADE singularities}\label{TabADE}
\end{table} 
In Table~\ref{TabADE} we indicate for each ADE singularity a non-degenerate invertible polynomial $f(x,y,z)$ with the correct Dolgachev and Gabrielov numbers. The corresponding polynomials  $f$ are self-dual and the Dolgachev and Gabrielov numbers of $f$ coincide. Moreover, the surface singularity given by $f=0$ has a minimal resolution with an exceptional configuration $\mathcal E$ whose dual graph is given by Fig.~\ref{FigExc}. Here $E$ and $E^i_1, \ldots, E^i_{\alpha_i-1}$ for $i=1,2,3$ are smooth rational curves of self-intersection number $-2$ and $\alpha_1, \alpha_2, \alpha_3$ are the Dolgachev numbers of $f$. 
In this case this graph coincides with the Coxeter-Dynkin diagram corresponding to a distinguished basis of vanishing cycles of this singularity. The graph is a classical Coxeter-Dynkin diagram of a root system of type $A_\mu$, $D_\mu$, $E_6$, $E_7$ or $E_8$. This is the reason why these singularities are called the ADE singularities. 

These singularities have many characterizations. 

The ADE singularities are just the quotients of $\CC^2$ by a finite subgroup $\Gamma \subset {\rm SL}(2;\CC)$. Let $\rho_0, \rho_1, \ldots, \rho_\ell$ be the equivalence classes of irreducible finite dimensional complex representations of $\Gamma$ where $\rho_0$ is the class of the trivial representation. J.~McKay \cite{McKay80} has observed that if $\rho: \Gamma \to {\rm SL}(2;\CC)$ is the given 2-dimensional representation of $\Gamma$ then the $(\ell+1) \times (\ell+1)$-matrix $B=(b_{ij})$, defined by decomposing the tensor products $\rho_j \otimes \rho= \bigoplus_i b_{ij}\rho_i$ into irreducible components, satisfies $B=2I-C$ where $C$ is the affine Cartan matrix of the corresponding root system. The Coxeter-Dynkin diagram corresponding to $C$ is just the extended Coxeter-Dynkin diagram of the corresponding root system. This is obtained by joining an additional vertex (corresponding to the trivial representation $\rho_0$) to the vertices $E$ and $E_1^3$ in the case $A_\mu$ ($(\alpha_1,\alpha_2, \alpha_3)=(1,1,\mu)$), to $E_2^3$ in the case $D_\mu$, $E_1^1$ in the case $E_6$, $E_1^2$ in the case $E_7$ and $E_1^3$ in the case $E_8$. It is shown in \cite{Ebeling96} that the corresponding diagram can be transformed to a diagram of type $T_{\alpha_1,\alpha_2, \alpha_3}$ by Gabrielov transformations (see Fig.~\ref{FigTpqr}). 

G.~Gonza\-lez-Sprinberg and J.-L.~Verdier \cite{GSV} and independently H.~Kn\"orrer \cite{Kn} gave a geometric interpretation of the McKay correspondence by identifying the Grothendieck group of the category of coherent sheaves on the minimal resolution with the representation ring of $\Gamma$. M.~Kapranov and E.~Vasserot \cite{KV} extended these results to the derived category of coherent sheaves on the minimal resolution, not just the Grothendieck group.

Let $A_f$ be the coordinate ring of the weighted homogeneous polynomial $f$. It is graded according to the weights of the variables. Let $P_f(t)$ be the Poincar\'e series of  the graded algebra $A_f$.
Let $c_{f,-}$ and $c_{f,0}$ be the Coxeter element of the root system and the affine root system associated to the singularity $f$ and $\phi_{f,-}(t)$ and $\phi_{f,0}(t)$ respectively the corresponding characteristic polynomial. The author has proved \cite{Ebeling02}:

\begin{theorem} \label{thm:Klein}
For a simple singularity not of type $A_{2k}$ we have
\[ P_f(t)= \frac{\phi_{f,-}(t)}{\phi_{f,0}(t)}. \]
\end{theorem}

R.~Stekolshchik has generalized this theorem to the root systems with non simply laced Coxeter-Dynkin diagrams \cite{Stek0, Stek1}.

The minimal resolution of an ADE singularity can be compactified to a rational surface $S_f$ containing the exceptional configuration $\mathcal E$. The author and D.~Ploog \cite{EP1} have given the following geometric realization of the graph $T_{\alpha_1, \alpha_2, \alpha_3}$. 

Let ${\rm Coh}(S_f)$ be the abelian category of coherent sheaves on $S_f$ and $K(S_f)$ its Grothendieck K-group. There is a natural bilinear pairing on $K(S_f)$ given by the Euler form  $\chi(A,B)=\sum_i(-1)^i\dim\Ext^i_{S_f}(A,B)$ for two coherent sheaves $A$ and $B$ 
on $S_f$. Let $N(S_f)$ be the numerical K-group which is obtained from $K(S_f)$ by dividing out the radical of the Euler form. Denote by ${\rm Coh}_{\mathcal E}(S_f)$ the abelian subcategory of ${\rm Coh}(S_f)$ consisting of sheaves whose support is contained in $\mathcal E$ and let $K_{\mathcal E}(S_f)$ be its K-group. The curves $E$ and $E^i_1, \ldots, E^i_{\alpha_i-1}$ for $i=1,2,3$ correspond to spherical objects in the category ${\rm Coh}_{\mathcal E}(S_f)$. (Recall that a coherent sheaf $F$ on $S_f$ is called {\em spherical} if
\[ \Ext^l_{S_f}(F,F) = \left\{ \begin{array}{ll}
                                    \CC, & l=0 \text{ or } l=2 \\
                                    0    & \text{else}
                           \end{array} \right.
   \qquad \text{ and } \qquad
   F\otimes\omega_{S_f}\cong F .)\]
It follows from \cite[Lemma~1]{EP1} that the Euler pairing between the classes 
\begin{equation}
[\calO_E(-1)],  [\calO_{E^1_1}(-1)], \ldots,  [\calO_{E^3_{\alpha_3-1}}(-1)], [\calO_E] \label{eq:spher1}
 \end{equation}
in $N(S_f)$ is encoded by the graph $T_{\alpha_1,\alpha_2, \alpha_3}$ (see Fig.~\ref{FigTpqr}) with the length of arms being equal to the Dolgachev numbers $\alpha_1, \alpha_2, \alpha_3$ of $f$. Using this description, the author and Ploog have given another proof of Theorem~\ref{thm:Klein} \cite{EP1}.

 In \cite{KST1}, H.~Kajiura, K.~Saito and A.~Takahashi proved the existence of a full strongly exceptional sequence in $D_{\rm Sg}^{L_f}(R_f)$ for a polynomial $f$ of ADE type. D.~Kussin, H.~Lenzing and H.~Meltzer \cite{KLM} discuss relations of these categories with weighted projective lines.

The unimodal singularities are the simple elliptic singularities $\widetilde{E}_6$,  $\widetilde{E}_7$ and $\widetilde{E}_8$ given by polynomials of type $T_{3,3,3}$, $T_{2,4,4}$ and $T_{2,3,6}$ respectively (where $a_f=0$), the singularities of type $T_{p,q,r}$ with $\Delta(p,q,r)>0$ and the 14 exceptional unimodal singularities (with $a_f=1$). Invertible polynomials for the simple elliptic singularities are given in \cite[Table~6]{ET1}. The singularities of type $T_{p,q,r}$ with $\Delta(p,q,r)>0$ are not weighted homogeneous. 

Now we come back to Arnold's strange duality. 

Let $f(x,y,z)$ be one of the invertible polynomials of Table~\ref{TabArnold}. Then $a_f=1$ and a Coxeter-Dynkin diagram of $f$ is given by the graph $S_{\gamma_1,\gamma_2,\gamma_3}$ which is an extension of the graph $T_{\gamma_1,\gamma_2,\gamma_3}$ by one vertex in accordance with Conjecture~\ref{Conj1}. Here $\gamma_1,\gamma_2,\gamma_3$ are the Gabrielov numbers of $f$.

H.~Pinkham \cite{Pinkham77} and  I.~V.~Dolgachev and V.~V.~Nikulin \cite{DN77} have shown that the Milnor fibre $V_f$ can be compactified in a weighted projective space $\PP(x,y,z,w)$ so that after minimal normal crossing resolution of singularities one obtains a K3-surface $S_f$. In this way, Arnold's strange duality can be considered as a special case of the mirror symmetry of K3-surfaces (see \cite{Dolgachev96}). One can use this K3-surface to find a categorical realization of a Coxeter-Dynkin diagram of the dual singularity, namely of the graph $S_{\alpha_1,\alpha_2,\alpha_3}$, where $\alpha_1,\alpha_2,\alpha_3$ are the Dolgachev numbers of $f$. This was obtained by the author and Ploog \cite{EP1}. 

Namely, the K3-surface $S_f$ carries an exceptional configuration $\mathcal E$ at $\infty$ whose dual graph is given by Fig.~\ref{FigExc}. Here $E$ and $E^i_1, \ldots, E^i_{\alpha_i-1}$ for $i=1,2,3$ are again smooth rational curves of self-intersection number $-2$ and  $\alpha_1, \alpha_2, \alpha_3$ are the Dolgachev numbers of $f$. 
The same construction as above can be applied to the K3 surface $S_f$. Moreover, the structure sheaf ${\mathcal O}_{S_f}$ of the K3-surface $S_f$ is also spherical. We consider in this case the classes
\begin{equation}
[\calO_E(-1)], [\calO_{E^1_1}(-1)], \ldots,  [\calO_{E^3_{\alpha_3-1}}(-1)],    [\calO_E], [{\mathcal O}_{S_f}]  \label{eq:spher2}
 \end{equation}
in $N(S_f)$. The Euler pairing between these classes is encoded by the graph $S_{\alpha_1,\alpha_2,\alpha_3}$ (see Fig.~\ref{FigSpqr}) with the lengths of arms being equal to the Dolgachev numbers $\alpha_1, \alpha_2, \alpha_3$ of $f$. 

Let $c_{f,+}$ be the Coxeter element corresponding to this graph and $\phi_{f,+}(t)$ be its characteristic polynomial. The author has shown \cite{Ebeling03}:
\begin{theorem} \label{thm:Fuchs}
\[ P_f(t)= \frac{\phi_{f,+}(t)}{\phi_{f,0}(t)}. \]
\end{theorem}
More generally, this result was proved for so-called Fuchsian singularities. These are normal surface singularities with a good $\CC^\ast$-action which are related to Fuchsian groups of the first kind. The hypersurface singularities among the Fuchsian singularities are just given by invertible polynomials $f$ with $a_f=1$. In this case, the orbifold curve  ${\mathcal C}_{(f,G_0)}$ is of the form $\HH/\Gamma$ where $\HH$ is the upper half plane and $\Gamma$ is a Fuchsian group of the first kind. The genus $g_{(f,G_0)}$ is called the genus of the Fuchsian singularity. There are 31 such singularities \cite{Sh,W}. There are 22 such singularities with genus $g_{(f,G_0)}=0$. They include the 14 exceptional unimodal singularities, the 6 heads of the bimodal series (see below) and two more. 

A possible generalization of the McKay correspondence for Fuchsian groups of genus 0 has been discussed by I.~Dolgachev \cite{Dolgachev09}.
 H.~Lenzing and J.~A.~de la Pe\~{n}a \cite{LP} derived Theorem~\ref{thm:Fuchs} for Fuchsian singularities of genus 0 from the representation theory of certain algebras related with weighted projective lines. In \cite{EP1, EP2} the author and Ploog derive this result for smoothable Fuchsian singularities of any genus from a generalization of the categorical approach indicated above.

In \cite{EP1}, the Coxeter elements were described as autoequivalences of triangulated categories as follows.  Denote by $D^b(S_f)$ the bounded derived category of coherent sheaves on $S_f$, 
and by ${\mathcal D}_{f,0} := {\mathcal D}^b_{\mathcal E}(S_f)$ the subcategory consisting of complexes whose support (of all homology sheaves) is contained in $\mathcal E$. These are triangulated categories and $D^b_{\mathcal E}(S_f)$ is a 2-Calabi-Yau category. We also consider the smallest full triangulated subcategory ${\mathcal D}_{f,+}$ of ${\mathcal D}^b(S_f)$ containing ${\mathcal D}_{f,0}$ and the structure sheaf ${\mathcal O}_{S_f}$  of $S_f$. Then the classes (\ref{eq:spher2}) form a spherical collection in ${\mathcal D}_{f,+}$. A spherical sheaf $F$ determines a {\em spherical twist} $\TTT_F: D^b(S_f) \to D^b(S_f)$ which is an autoequivalence of the category \cite{ST}. Consider the autoequivalences
\begin{eqnarray*}
c_{f,0} & := & \TTT_\tytt{\calO_E(-1)} \TTT_\tytt{\calO_{E_1^{1}}(-1)} \cdots \TTT_\tytt{\calO_{E_{\alpha_r-1}^{r}}(-1)}  \TTT_\tytt{\calO_E}, \\
c_{f,+}     & := & c_{f,0}               \TTT_\tytt{\calO_{S_f}}.
\end{eqnarray*}
The autoequivalences $c_{f,0}$ and $c_{f,+}$ correspond to the Coxeter elements of the graphs $T_{\alpha_1,\alpha_2,\alpha_3}$ and $S_{\alpha_1,\alpha_2,\alpha_3}$ respectively. 

In \cite{KST2}, H.~Kajiura, K.~Saito and A.~Takahashi proved the existence of a full strongly exceptional sequence in $D_{\rm Sg}^{L_f}(R_f)$ for a weighted homogeneous polynomial $f$  of with $\eps_f=-a_f=-1$ and genus $g_{(f,G_0)}=0$. This includes the case of the 14 exceptional unimodal singularities. Lenzing and de la Pe\~{n}a \cite{LP} proved that the category $D_{\rm Sg}^{L_f}(R_f)$ in this case is equivalent to the derived category of finitely generated modules over the extended canonical algebra associated with the weighted projective line ${\mathcal C}_{(f,G_f)}$. The relation between the categories $D_{\rm Sg}^{L_f}(R_f)$ and ${\mathcal D}_{f,+}$  for the 14 exceptional unimodal singularities was studied by M.~Kobayashi, M.~Mase and K.~Ueda in \cite{KMU}.

We now turn to the bimodal singularities. They are also classified by Arnold \cite{Ar}. They fall into the following 8 infinite series of singularities (for $k \in \ZZ$)
\[ J_{3,k}, Z_{1,k}, Q_{2,k}, W_{1,k},  S_{1,k}, U_{1,k}, \  k \geq 0, \mbox{ and } W^\sharp_{1,k}, S^\sharp_{1,k}, \ k \geq 1,
\]
and the 14 exceptional singularities
\[ E_{18}, E_{19}, E_{20}, Z_{17}, Z_{18}, Z_{19}, Q_{16}, Q_{17}, Q_{18}, W_{17}, W_{18}, S_{16}, S_{17}, U_{16}.
\]
One can find weighted homogeneous polynomials for the classes for $k=0$ in the series and for the 14 exceptional singularities. In each of these classes, one can find a non-degenerate invertible polynomial $f$. These polynomials, their Berglund-H\"ubsch transposes and their Dolgachev numbers $A_{(f, G_f)}=(\alpha_1,\alpha_2, \alpha_3)$ and Gabrie\-lov numbers $\Gamma_{(f, \{ e \})}=(\gamma_1, \gamma_2, \gamma_3)$ are indicated in Table~\ref{TabBi}. Note that the dual singularities are only bimodal for the self-dual cases, in the other cases other singularities are involved.
\begin{table}[h]
\begin{center}
\begin{tabular}{cccccc}
\hline
Name & $\alpha_1,\alpha_2, \alpha_3$  & $f$   & $\widetilde{f}$ & $\gamma_1, \gamma_2, \gamma_3$ & Dual \\
\hline
$J_{3,0}$ & $2,4,6$ & $x^6y+y^3+z^2$& $x^6+xy^3+z^2$  & $2,3,10$ & $Z_{13}$ \\
$Z_{1,0}$ & $2,4,8$ & $x^5y + xy^3 +z^2$ & $x^5y + xy^3 +z^2$ & $2,4,8$ & $Z_{1,0}$ \\
$Q_{2,0}$ & $2,4,10$ & $x^4y + y^3 + xz^2$ & $x^4z + xy^3 +z^2$ & $3,3,7$ & $Z_{17}$\\
$W_{1,0}$ & $2,6,6$  & $x^6+y^2+yz^2$ & $x^6+y^2z+z^2$ & $2,6,6$ & $W_{1,0}$ \\
$S_{1,0}$ & $2,6,8$ & $x^5+xy^2+yz^2$ & $x^5y+y^2z+z^2$ & $3,5,5$ & $W_{17}$ \\
$U_{1,0}$ & $3,4,6$  & $x^3+xy^2+yz^3$ & $x^3y+y^2z+z^3$ & $3,4,6$ & $U_{1,0}$ \\
\hline
$E_{18}$ & $3, 3, 5$  & $x^5z + y^3 + z^2$ & $x^5 + y^3 + xz^2$ & $2,3,12$ & $Q_{12}$ \\
$E_{19}$ & $2, 4, 7$  & $x^7y + y^3+z^2$  & $x^7 + xy^3+z^2$ & $2,3,12$ & $Z_{1,0}$  \\
$E_{20}$ & $2, 3, 11$ & $x^{11} +y^3 +z^2$ & $x^{11} +y^3 +z^2$ & $2, 3, 11$ & $E_{20}$ \\
$Z_{17}$ & $3, 3, 7$ & $x^4z + xy^3 + z^2$ & $x^4y + y^3 + xz^2$ & $2,4,10$ & $Q_{2,0}$  \\
$Z_{18}$ & $2, 4, 10$  & $x^6y + xy^3 + z^2$ & $x^6y + xy^3 + z^2$ & $2, 4, 10$ & $Z_{18}$  \\
$Z_{19}$ & $2, 3, 16$ & $x^9 + xy^3 + z^2$  & $x^9y + y^3 + z^2$ & $2, 4, 9$ & $E_{25}$\\
$Q_{16}$ & $3, 3, 9$ & $x^4z + y^3 + xz^2$ & $x^4z + y^3 + xz^2$ & $3, 3, 9$ & $Q_{16}$ \\
$Q_{17}$ & $2, 4, 13$ & $x^5y+y^3+xz^2$  & $x^5z+xy^3+z^2$ & $3,3,9$  & $Z_{2,0}$  \\
$Q_{18}$ & $2, 3, 21$ & $x^8+y^3+xz^2$ & $x^8z+y^3+z^2$ & $3, 3, 8$ & $E_{30}$  \\
$W_{17}$ & $3, 5, 5$ & $x^5z+yz^2+y^2$ & $x^5+xz^2+y^2z$ & $2,6,8$ & $S_{1,0}$  \\
$W_{18}$ & $2, 7, 7$ & $x^7+ y^2+ yz^2$ & $x^7+ y^2z+ z^2$ & $2, 7, 7$ & $W_{18}$  \\
$S_{16}$ & $3, 5, 7$ & $x^4y+xz^2+y^2z$ & $x^4y+xz^2+y^2z$ & $3, 5, 7$ & $S_{16}$  \\
$S_{17}$ & $2, 7, 10$  & $x^6+xy^2+yz^2$ & $x^6y+y^2z+z^2$ & $3,6,6$ & $X_{2,0}$  \\
$U_{16}$ & $5, 5, 5$ & $x^5+y^2z+yz^2$ & $x^5+y^2z+yz^2$ & $5, 5, 5$ & $U_{16}$  \\
\hline
\end{tabular}
\end{center}
\caption{Strange duality of the bimodal singularities} \label{TabBi}
\end{table}

The Gorenstein parameter $a_f$ takes the value 1 for the classes for $k=0$ in the series and it takes the values 2,3 and 5 for the exceptional bimodal singularities. 
Coxeter-Dynkin diagrams with respect to distinguished bases of vanishing cycles for the bimodal singularities have been computed in \cite{Ebeling83}. Each of these diagrams is an extension of the corresponding graph  $T_{\gamma_1,\gamma_2,\gamma_3}$ by $a_f$ vertices in accordance with Conjecture~\ref{Conj1}.
The author and Ploog \cite{EP3} have given a geometric construction of these Coxeter-Dynkin diagrams by a procedure as above using compactifications of the Milnor fibres of the polynomials $\widetilde{f}$. (Note that there are some mistakes in \cite[Table~4]{EP3}, they are corrected in arXiv:1102.5024.) M.~Mase and K.~Ueda \cite{MU} have shown that this result can also be derived from Conjecture~\ref{Conj1} and they used our construction to relate the Berglund-H\"ubsch duality in this case to Batyrev's polar duality for the associated toric K3 surfaces.

We also mention some other examples given in \cite{ET2}.
Consider the invertible polynomial $f(x,y,z)=x^2+xy^3+yz^5$. The canonical system of weights is $W_f=(15,5,5;30)$, so $c_f=5$ and the reduced system of weights is $W_f^{\rm red}=(3,1,1;6)$. This is again a singularity with $a_f=1$, but the genus of $\calC_{(f, G_0)}$ is equal to two and there are no isotropic points. The Dolgachev numbers of the pair $(f,G_f^{\rm fin})$ are $(5,5,5)$ and $G_0^T$ is generated by the element $({\bf e}[\frac{1}{5}], {\bf e}[\frac{3}{5}], {\bf e}[\frac{1}{5}])$. This group has two elements of age 1. The singularity $f^T(x,y,z)-xyz$ is right equivalent to the cusp singularity $x^5+y^5+z^5-xyz$. We recover the example of Seidel \cite{Se5}. Similarly, let $f(x,y,z)=x^3y+y^3z+z^3x$. Then again $a_f=1$, the genus of  $C_{(f, G_0)}$ is equal to three and there are no isotropic points. The group $G_0^T$ is generated by the element $({\bf e}[\frac{1}{7}], {\bf e}[\frac{2}{7}], {\bf e}[\frac{4}{7}])$ and $f^T(x,y,z)-xyz$ is right equivalent to the cusp singularity $x^7+y^7+z^7-xyz$.

More generally, let $g$ be an integer with $g \geq 2$ and consider the invertible polynomial $f(x,y,z)=x^{2g+1}+y^{2g+1}+z^{2g+1}$ together with the group $G$ generated by $({\bf e}[\frac{1}{2g+1}], {\bf e}[\frac{1}{2g+1}], {\bf e}[\frac{2g-1}{2g+1}])$. Then the genus of the curve $C_{(f,G)}$ is equal to $g$ and we recover the examples of A.~Efimov \cite{Efimov:2009}.

\section{Complete intersection singularities as mirrors} \label{sect:EW}
We shall now derive the extension of Arnold's strange duality discovered by the author and C.~T.~C.~Wall  from the mirror symmetry and the Berglund-H\"ubsch transposition of invertible polynomials. This is the contents of the paper \cite{ET5}.

The invertible polynomials $f(x,y,z)$ given in Table~\ref{TabBi} for the singularities  $J_{3,0}$, $Z_{1,0}$, $Q_{2,0}$, $W_{1,0}$, $S_{1,0}$ and $U_{1,0}$, which are the heads of the bimodal series, satisfy $[G_f : G_0]=2$. There corresponds an action of $\widetilde{G}_0=\ZZ/2\ZZ$ on the transpose polynomial $\widetilde{f}$. We choose the coordinates such that this action is given by $(x,y,z) \mapsto (-x,-y,z)$. The Dolgachev numbers of the pairs $(f, G_0)$ are given in Table~\ref{TabBi0} (see also \cite[Table~4]{ET2}). Moreover, there is a one-parameter family $F$ (depending on a complex parameter $a$) of weighted homogeneous polynomials defining these singularities, it is also indicated in Table~\ref{TabBi0}. It is natural from the mirror symmetry view point to expect that adding one monomial to an invertible polynomial is dual to having another $\CC^\ast$-action on the dual polynomial. This will be elaborated in the sequel.
\begin{table}[h]
\begin{center}
\begin{tabular}{ccc}
\hline
Name   & $F$   & $A_{(f,G_0)}$  \\
\hline
$J_{3,0}$ & $x^3+xy^6+z^2+ax^2y^3$, $a \neq \pm 2$ & $2,2,2,3$ \\
$Z_{1,0}$  & $x^5y + xy^3 +z^2+ax^2y^3$, $a \neq \pm 2$ & $2,2,2,4$ \\
$Q_{2,0}$ & $x^3 + xy^4 + yz^2+ax^2y^2$, $a \neq \pm 2$ & $2,2,2,5$ \\
$W_{1,0}$ & $x^6+y^2+yz^2+ax^3y$, $a \neq \pm 2$ & $2,2,3,3$  \\
$S_{1,0}$  & $x^5+xy^2+yz^2+ax^3y$, $a \neq \pm 2$ & $2,2,3,4$ \\
$U_{1,0}$  & $x^3+xy^2+yz^3+ax^2y$, $a \neq \pm 2$ & $2,3,3,3$ \\
\hline
\end{tabular}
\end{center}
\caption{Heads ($k=0$) of bimodal series} \label{TabBi0}
\end{table}

The non-degenerate invertible polynomials $f(x,y,z)$ with $[G_f : G_0]=2$ are classified in \cite[Proposition~5]{ET5}. There are 5 possible types each depending on parameters $p_1,p_2,p_3$ or $p_1,q_2,q_3$ subject to certain conditions.

\begin{sloppypar}
Let $L_0$ be the quotient of $L_{f}$ corresponding to the subgroup $\widehat{G}_0$ of $\widehat{G}_{f}$ (cf.\ Sect~\ref{sect:inv}). 
We consider $4 \times 3$-matrices $E = (E_{ij})^{i=1,2,3,4}_{j=1,2,3}$ such that 
\[ \ZZ \vec{x} \oplus \ZZ \vec{y} \oplus \ZZ \vec{z} \oplus \ZZ \vec{f} / \langle E_{i1} \vec{x} +E_{i2} \vec{y}+E_{i3} \vec{z}= \vec{f}, i=1, \ldots ,4 \rangle \cong L_0
\]
and $\calC_{(F,G_0)} := [ (F^{-1}(0) \setminus \{ 0 \})/ \widehat{G}_0]$, where $F:= \sum_{i=1}^4 a_i x^{E_{i1}}y^{E_{i2}}z^{E_{i3}}$, is a smooth projective line with 4 isotropic points whose orders are $\alpha_1, \alpha_2, \alpha_3, \alpha_4$, where $A_{(f,G_0)}=(\alpha_1, \alpha_2, \alpha_3, \alpha_4)$ are the Dolgachev numbers of the pair $(f,G_0)$ defined above, for general $a_1, a_2, a_3,  a_4$. These matrices are classified in \cite[Proposition~2]{ET5}.
\end{sloppypar}

We associate to these matrices a pair of polynomials as follows. We observe that the kernel of the matrix $E^T$ is either generated by the vector $(1,1,0,-2)^T$ or by the vector $(1,1,-1,-1)^T$. 
Let $R:=\CC[x,y,z,w]$. In the first case, there exists a $\ZZ$-graded structure on $R$ given by the $\CC^\ast$-action
\[ \lambda \ast (x,y,z,w) = (\lambda x, \lambda y, z , \lambda^{-2} w) \quad \mbox{for } \lambda \in \CC^\ast.
\]
In the second case, there exists a $\ZZ$-graded structure on $R$ given by the $\CC^\ast$-action
\[ \lambda \ast (x,y,z,w) = (\lambda x, \lambda y, \lambda^{-1}z , \lambda^{-1} w) \quad \mbox{for } \lambda \in \CC^\ast.
\]
Let $R= \bigoplus_{i \in \ZZ} R_i$ be the decomposition of $R$ according to one of these $\ZZ$-gradings. Let $E^T$ be the transposed matrix.  We associate to this the polynomial 
\[\widetilde{f}(x,y,z,w) := x^{E_{11}}y^{E_{21}}z^{E_{31}}w^{E_{41}}+ x^{E_{12}}y^{E_{22}}z^{E_{32}}w^{E_{42}} +  x^{E_{13}}y^{E_{23}}z^{E_{33}}w^{E_{43}}.
\]
In the first case, we have $\widetilde{f} \in R_0=\CC[x^2w,y^2w,z,xyw]$.
Let
\[ X:=x^2w, \quad Y:=y^2w, \quad Z:=z, \quad W:=xyw. 
\]
In these new coordinates, we obtain a pair of polynomials
\[
\widetilde{\ff}_1(X,Y,Z,W) = XY-W^2, \quad \widetilde{\ff}_2(X,Y,Z,W) = \widetilde{f}(X,Y,Z,W).
\]

In the second case, we have  $\widetilde{f} \in R_0=\CC[xw,yz,xz,yw]$. Let
\[ X:=xw, \quad Y:=yz, \quad Z:=xz \quad W:=yw. 
\]
In these new coordinates, we obtain a pair of polynomials
\[
\widetilde{\ff}_1(X,Y,Z,W) = XY-ZW, \quad \widetilde{\ff}_2(X,Y,Z,W) = \widetilde{f}(X,Y,Z,W).
\]

Now we choose for each of the matrices $E$ special values $a_1, a_2, a_3, a_4$ such that the corresponding polynomial $F$ has a non-isolated singularity. We denote this polynomial by $\ff$. In two cases, we already have additional conditions on the parameters. In the remaining cases, we consider conditions on the parameters $p_1,p_2,p_3$ ($p_1,q_2,q_3$) such that the polynomial $\ff(x,y,z)$ is of the form
\[ 
\ff(x,y,z)= u(x,y,z)+v(x,y,z)(x-y^e)^2
\]
or
\[ 
\ff(x,y,z)= u(x,y,z)+v(x,y,z)(y-x^e)^2
\]
for some monomials $u(x,y,z)$ and $v(x,y,z)$ and some integer $e \geq 2$. We
consider the cusp singularity $\ff(x,y,z)-xyz$ and perform the coordinate change $x \mapsto x+y^e$ or $y \mapsto y+x^e$ respectively. Then $\ff(x,y,z)-xyz$ is transformed to $\hh(x,y,z)-xyz$. Some of the new polynomials $\hh$ have 4 monomials and others only 3. We restrict our consideration to the cases where the polynomial $\hh$ has 4 monomials.
The singularities defined by the polynomials $\hh(x,y,z)$ will be called {\em virtual singularities}. 

We summarize our duality in Table~\ref{TabVCIS}.
\begin{table}[h]
\begin{center}
\begin{tabular}{ccc}
\hline
Type   & $\hh$ & $(\widetilde{\ff}_1, \widetilde{\ff}_2)$\\
\hline
$\begin{array}{c} {\rm IIA} \\ p_2=3 \end{array}$ &   $-y^{\frac{p_3}{6}+1}z  +z^{p_1} +x^3+x^2y^{\frac{p_3}{6}}$ & 
$\left\{ \begin{array}{c} XY-W^2\\ XW+ Y^{\frac{p_3}{6}}+Z^{p_1}\end{array}\right\}$\\
$\begin{array}{c} {\rm IIB} \\ p_2=2 \end{array}$ &  $-x^{\frac{p_1}{2}+1}z+ y^2 + yz^{\frac{p_3}{2}} +x^{\frac{p_1}{2}}z^{\frac{p_3}{2}}$ & $\left\{ \begin{array}{c} XY-W^2\\ X^{\frac{p_1}{2}}+YZ+Z^{\frac{p_3}{2}} \end{array}\right\}$\\
$\begin{array}{c} {\rm IIB}^\sharp \\ p_2=2 \end{array}$ &  $-x^{\frac{p_1}{2}+1}z+y^2 + yz^{\frac{p_3}{2}} +x^{\frac{p_1}{2}}y$ & $\left\{ \begin{array}{c} XY-ZW\\ X^{\frac{p_1}{2}}+YW+Z^{\frac{p_3}{2}} \end{array}\right\}$\\
$\begin{array}{c} {\rm III} \\ q_2=2 \end{array}$ &   $-y^{\frac{q_3}{2}+1}z +z^{p_1}  + x^3y +x^2y^{\frac{q_3}{2}+1}$ & $\left\{ \begin{array}{c} XY-W^2\\ (X+Y^{\frac{q_3}{2}})W+Z^{p_1} \end{array}\right\}$\\
$\begin{array}{c} {\rm IV}_1\\p_1=3 \end{array}$ &   $-y^{\frac{p_2}{6}+1}z+x^3+yz^{\frac{p_3}{p_2}}+x^2y^{\frac{p_2}{6}}$ & $\left\{ \begin{array}{c} XY-W^2\\ XW+Y^{\frac{p_2}{6}}Z+ Z^{\frac{p_3}{p_2}} \end{array}\right\}$\\
$\begin{array}{c} {\rm IV}_2\\ \frac{p_2}{p_1}=2 \end{array}$ &  $-x^{\frac{p_1+1}{2}}z+ xy^2 +yz^{\frac{p_3}{p_2}}+x^{\frac{p_1-1}{2}}z^{\frac{p_3}{p_2}}$ & $\left\{ \begin{array}{c} XY-W^2\\ X^{\frac{p_1-1}{2}}W+YZ+ Z^{\frac{p_3}{p_2}} \end{array}\right\}$\\
$\begin{array}{c} {\rm IV}_2^\sharp\\ \frac{p_2}{p_1}=2 \end{array}$ &  $-x^{\frac{p_1+1}{2}}z+ xy^2 +yz^{\frac{p_3}{p_2}}+x^{\frac{p_1+1}{2}}y$ & $\left\{ \begin{array}{c} XY-ZW\\ X^{\frac{p_1-1}{2}}W+YW+Z^{\frac{p_3}{p_2}}\end{array}\right\}$\\
\hline
\end{tabular}
\end{center}
\caption{Duality between virtual singularities and complete intersection singularities} \label{TabVCIS}
\end{table}

One can associate Dolgachev and Gabrielov numbers to the virtual singularities and the pairs $(\widetilde{\ff}_1, \widetilde{\ff}_2)$ in an analogous way as above, see \cite[Sect.~4]{ET5}. Here the Gabrielov numbers of the virtual singularities and the Dolgachev numbers of the pairs $(\widetilde{\ff}_1, \widetilde{\ff}_2)$ are triples, but the Dolgachev numbers of the virtual singularities and the Gabrielov numbers of the pairs $(\widetilde{\ff}_1, \widetilde{\ff}_2)$ are quadruples of numbers which are divided into two pairs in a natural way. One obtains the following theorem \cite[Theorem~4]{ET5}:
\begin{theorem}[E., Takahashi]  \label{thm:duality}
The Gabrielov numbers of the polynomial $\hh$ corresponding to a virtual singularity  coincide with the Dolgachev numbers of the dual pair $(\widetilde{\ff}_1,\widetilde{\ff}_2)$ and, vice versa, the Gabrielov numbers of a pair $(\widetilde{\ff}_1,\widetilde{\ff}_2)$ coincide with the Dolgachev numbers of the dual polynomial $\hh$.
\end{theorem}
There is also an extension of Saito's duality to this duality, see \cite[Corollary~6]{ET5}.

Now we consider again the cases with small Gorenstein parameter $a_f$: For the case $a_f <0$ see \cite[Table~8]{ET5}. There are no non-degenerate invertible polynomials with $[G_f : G_0]=2$ and $a_f=0$. 
Next consider the case $a_f=1$. It turns out that the virtual singularities in this case are exactly the virtual singularities corresponding to the bimodal series. Namely, by setting $k=-1$ in Arnold's equations one obtains polynomials which are similar to our polynomials $\hh$ for certain types and parameters $p_1,p_2,p_3$ or $p_1,q_2,q_3$ respectively. These types and parameters are listed in Table~\ref{TabBimon} together with the corresponding Dolgachev and Gabrielov numbers of $\hh$ and the corresponding dual pairs $(\widetilde{\ff}_1, \widetilde{\ff}_2)$ defining an isolated complete intersection singularity (ICIS).
\begin{table}[h]
\begin{center}
\begin{tabular}{ccccccc}
\hline
Name & Type & ${\rm Dol}(\hh)$ &  ${\rm Gab}(\hh)$  &  ${\rm Dol}(\widetilde{\ff}_1, \widetilde{\ff}_2)$ &  ${\rm Gab}(\widetilde{\ff}_1, \widetilde{\ff}_2)$ & Dual  \\
\hline
$J_{3,-1}$  & IIA, $2,3,18$ & $2,2;2,3$ & $2,3,10$  & $2,3,10$ & $2,2;2,3$  & $J_9'$  \\
$Z_{1,-1}$ & III, $2,2,4$ &  $2,2;2,4$  & $2,4,8$ & $2,4,8$ & $2,2;2,4$ & $J_{10}'$  \\
$Q_{2,-1}$ & ${\rm IV}_1$, $3,12,24$ & $2,2;2,5$  & $3,3,7$  & $3,3,7$ & $2,2;2,5$ & $J_{11}'$  \\
$W_{1,-1}$ & IIB, $6,2,4$ & $2,2;3,3$ & $2,6,6$  & $2,6,6$ & $2,2;3,3$ & $K_{10}'$ \\
$W^\sharp_{1,-1}$ & ${\rm IIB}^\sharp$, $6,2,4$ & $2,3;2,3$  & $2,5,7$ & $2,5,7$ & $2,3;2,3$ & $L_{10}$ \\
$S_{1,-1}$ & ${\rm IV}_2$, $5,10,20$  & $2,2;3,4$  & $3,5,5$  & $3,5,5$ & $2,2;3,4$ & $K_{11}'$ \\
$S^\sharp_{1,-1}$ & ${\rm IV}_2^\sharp$, $5,10,20$ & $2,3;2,4$  & $3,4,6$  & $3,4,6$ & $2,3;2,4$ & $L_{11}$ \\
$U_{1,-1}$ & ${\rm IV}_2^\sharp$, $3,6,18$ & $2,3;3,3$ & $4,4,5$  & $4,4,5$ & $2,3;3,3$ & $M_{11}$  \\
\hline
\end{tabular}
\end{center}
\caption{Strange duality between virtual bimodal singularities and ICIS} \label{TabBimon}
\end{table}

Let $\hh(x,y,z)=0$ be the equation for one of the virtual bimodal singularities. It turns out that $\hh$ has besides the origin an additional critical point which is of type $A_1$. One can find a Coxeter-Dynkin diagram with respect to a distinguished basis of vanishing cycles corresponding to all the critical points of the form $S_{\gamma_1,\gamma_2,\gamma_3}$ where $\gamma_1,\gamma_2,\gamma_3$ are the Gabrielov numbers  of $\hh$. 

To a graph of type $S_{\gamma_1,\gamma_2,\gamma_3}$, there corresponds an extended canonical algebra in the sense of \cite{LP}. The 14 cases of the exceptional unimodal singularities (see Table~\ref{TabArnold}) and the 8 cases of Table~\ref{TabBimon} correspond to those extended canonical algebras where  the number $t$ of weights is equal to 3 and the Coxeter element is semi-simple and has only roots of unity as eigenvalues (cf.\ \cite[Theorem~3.4.3 and Table~3.4.2]{Ebeling87} and \cite[Table~2]{LP}).
\section{Ackowledgements}
I would like to thank the referee for useful comments.

\section{References}

{}
\renewcommand{\refname}{}    
\vspace*{-20pt}              

\end{document}